\documentclass{conm-p-l}
\usepackage{graphicx,amssymb,bbm}
\usepackage{rotating}
\usepackage[mathscr]{eucal} 
\def\addforarchive{\begin{picture}(0,0)
                   \put(60,250){{\small\sf math.CT/0511590}}
                   \put(60,235){{\small\sf KCL-MTH-05-14}}
                   \put(60,220){{\small\sf Hamburger$\;$Beitr\"age$\;$zur$\;$Mathematik$\;$Nr.$\;$222}}
                   \put(60,205){{\small\sf ZMP-HH/05-06}}
                   \eP} 
\copyrightinfo{2005}{Ribbon Graphs United}

\newtheorem{theorem}{Theorem}
\newtheorem{lemma}[theorem]{Lemma}
\newtheorem{proposition}[theorem]{Proposition}
\newtheorem{corollary}[theorem]{Corollary}
\theoremstyle{definition}
\newtheorem{definition}[theorem]{Definition}
\theoremstyle{remark}
\newtheorem{remark}[theorem]{Remark}

\def\AA            {{\!A}}
\newcommand\aAa[2] {{}_{#1}{A}_{#2}}
\def\alg           {algebra}
\def\alinm         {\alpha^-_\AA}
\def\alinp         {\alpha^+_\AA}
\def\alinpm        {\alpha^\pm_\AA}
\def\Aop           {\ensuremath{A^{\rm opp}_{}}}
\def\atimes        {\mbox{\small$\#$}}
\def\ati           {\,{\atimes}\,}
\def\Atop          {\ensuremath{A_\topsy}}
\def\Atopx         {\ensuremath{A_\topsy^\times}}
\def\Aut           {{\rm Aut}}
\def\auto          {automorphism}
\def\be            {\begin{equation}}
\def\bearll        {\begin{array}{ll}}
\def\bPo           {\begin{picture}(0,0)}
\def\Br            {{\rm FBr}}
\def\C             {\ensuremath{\mathcal C}}
\def\CA            {\ensuremath{\mathcal C_{\!A|}}}
\def\CAA           {\ensuremath{\mathcal C_{\!A|A}}}
\def\CAB           {\ensuremath{\mathcal C_{\!A|B}}}
\def\CAr           {\ensuremath{\mathcal C_{|A}}}
\def\CB            {\ensuremath{\mathcal C_{B|}}}
\def\cala          {\ensuremath{\mathcal A}}

\def\calx          {\ensuremath{\mathcal X}}
\def\calxc         {\ensuremath{\mathcal X_{\CC}}}
\def\calyc         {\ensuremath{\mathcal Y_{\C}}}
\def\cat           {category}
\def\cats          {categories}
\def\C             {\ensuremath{\mathcal C}}
\def\CC            {\ensuremath{\mathscr O}}
\def\Cf            {\mbox{\sl Cor}}
\def\CfA           {\ensuremath{\Cf_{\!A}}}
\def\cft           {conformal field theory}
\def\cfta          {\ensuremath{\mbox{\sl cft}_{\C,A}}}
\def\cftf          {\ensuremath{\mbox{\sl cft}_{\CC}}}
\def\cftco         {\ensuremath{\mbox{\sl c-cft}_{\C}}}
\def\cfts          {conformal field theories}
\def\cir           {\,{\circ}\,}
\def\complex       {\ensuremath{\mathbbm C}}
\def\complexx      {{\ensuremath{{\mathbbm C}^{\times}_{}}}}
\def\corfu         {correlation function}
\def\Corfu         {Correlation function}
\def\Cst           {\ensuremath{\mathcal C^*_{}}}

\def\Ctop          {\ensuremath{C_\topsy}}
\def\Ctopx         {\ensuremath{C_\topsy^\times}}
\newcommand\dictentre[2] {#1 && #2 \\}
\newcommand\dictentry[2] {#1 && #2 \\[.3em]}
\def\dim           {{\rm dim}}
\def\dimc          {{\rm dim}_\complex}
\def\dimk          {{\rm dim}_\koerper}
\def\ee            {\end{equation}}
\def\eear          {\end{array}}
\def\End           {\ensuremath{\mathrm{End}}} 
\def\eP            {\end{picture}}
\def\eps           {\varepsilon}
\def\epsa          {\epsilon_\alpha}
\def\epsap         {\epsilon_{\alpha'}}
\def\Epsilon       {\nu}
\def\epsnat        {\varepsilon_\natural}
\def\eq            {\,{=}\,}
\newcommand\erf[1] {(\ref{#1})}
\def\FBg           {Fro\-be\-ni\-us\hy Brau\-er group}
\def\fsi           {Fro\-be\-ni\-us\hy Schur indicator}
\def\Fv            {\ensuremath{F_\nu}}   
\def\Fx            {\ensuremath{F_\chi}}  
\def\Hom           {\ensuremath{\mathrm{Hom}}} 
\def\HomA          {\ensuremath{\mathrm{Hom}_{\!A}}}
\def\HomAA         {\ensuremath{\mathrm{Hom}_{\!A|A}}}
\def\HomAB         {\ensuremath{\mathrm{Hom}_{\!A|B}}}
\def\hy            {$\mbox{-\hspace{-.66 mm}-}$\linebreak[0]}
\def\I             {{\mathcal I}}
\def\id            {\mbox{\sl id}}
\def\idA           {\id_A}
\def\idss          {\mbox{\scriptsize\sl id}}
\def\Im            {{\rm Im}\,}
\def\iN            {\,{\in}\,}
\def\Inn           {{\rm Inn}}
\newcommand\inclpic[3] {{\begin{picture}{#3} \scalebox{.#2}
                   {\includegraphics{sydL_#1.eps}}
                   \setlength{\unitlength}{1pt} \eP}}
\newcommand\inn[1] {\varpi_{#1}}
\def\intro         {introduction }
\def\io            {\mbox{\sl\i}}
\def\iox           {\ensuremath{\io(\X)}}
\def\koerper       {\ensuremath{\Bbbk}} 
\def\M             {\ensuremath{\mathrm M}}
\def\ML            {\ensuremath{\mathcal M}}
\def\MR            {\ensuremath{\mathcal M'}}
\def\mtc           {modular tensor category}
\def\mtcs          {modular tensor categories}
\def\multo         {*}
\def\MX            {\ensuremath{{\mathrm M}_{\mathrm X}}}
\def\myitems       {\def\leftmargini{2.2em}~\\[-1.5em]
                   \begin{itemize}\addtolength\itemsep{2pt}} 
\newcommand\N[3]   {{N_{#1#2}}^{\!\!#3}}
\def\obj           {{\mathcal O}bj}
\def\one           {{\bf1}}
\def\ort           {\ensuremath{\mathrm{or}_2}}
\def\ota           {\,{\otimes}_{\!\!A}^{}\,}
\def\otA           {\ota}
\def\Ota           {{\otimes}_{\!A}^{}}
\def\oti           {\,{\otimes}\,}
\def\Oti           {{\otimes}}
\def\Otim          {{\otimes^{\!-}}}
\def\Otip          {{\otimes^{\!+}}}
\def\parfu         {partition function}
\def\PicaC         {\ensuremath{\mathcal{P}\!ic(\C)}}
\def\PicC          {\ensuremath{\mathrm{Pic}(\C)}}
\def\PicCAA        {\ensuremath{\mathrm{Pic}(\CAA)}}
\def\PICCAA        {\ensuremath{\mathcal P\!\mbox{\sl ic}(\CAA)}}
\def\PicoC         {\ensuremath{\mathrm{Pic}^\circ(\C)}}
\def\q             {quantum }
\def\Q             {Quantum }
\def\qft           {quantum field theory}
\def\qfts          {quantum field theories}
\def\reals         {{\mathbb R}}

\def\repV          {\ensuremath{\mathcal{R}ep(\V)}}

\def\rep           {representation}

\newcommand\setulen[2]{\setlength\unitlength{.#1#2pt}}
\def\simj          {\backsim}  
\def\ssfa          {symmetric special Fro\-be\-ni\-us algebra}
\def\ssFA          {\ssfa}
\def\sss           {\scriptscriptstyle}
\def\tft           {topological field theory}
\def\tftc          {\ensuremath{\mbox{\sl tft}_{\C}}}
\def\tfts          {topological field theories}
\def\Times         {\,{\times}\,}
\def\To            {\,{\to}\,}
\def\topsy         {{\circ}}
\newcommand\turlabl[1]{{\begin{turn}{90}$\sss #1$\end{turn}}}
\newcommand\Turlabl[1]{{\begin{turn}{270}$\sss #1$\end{turn}}}
\def\V             {\ensuremath{\mathscr V}}
\def\Vect          {\ensuremath{\mathcal Vect_\complex}}
\def\Vectp         {\ensuremath{\mathcal Vect_\complex^{\sss\bullet}}}
\def\verylongrarr  {{{-}\!\!{-}\!\!\!\!\longrightarrow}}
\def\vPhi          {\varPhi}
\def\vPsi          {\varPsi}
\def\vTheta        {\varTheta}
\def\X             {\ensuremath{\mathrm X}}
\def\Xh            {\ensuremath{\widehat{\mathrm X}}}

\def\zet           {\ensuremath{\mathbb Z}}
\def\Zu            {{\rm Z}}

\begin{document}

\title[Ribbon categories, CFT and Frobenius algebras]
 {Ribbon categories and (unoriented) CFT:\\
 Frobenius algebras, automorphisms, reversions}

\author[J.~Fuchs]{J\"urgen Fuchs}
\address{Institutionen f\"or fysik, Karlstads Universitet,
Universitetsg.~5, S--65188 Karlstad}
\email{jfuchs@fuchs.tekn.kau.se}

\author[I.~Runkel]{Ingo Runkel}
\address{Department of Mathematics, King's College London,
Strand, London WC2R 2LS}
\email{ingo@mth.kcl.ac.uk}

\author[C.~Schweigert]{Christoph Schweigert}
\address{Fachbereich Mathematik, Universit\"at Hamburg,
Bundesstra{\rm\ss{}}e 55, D--20146 Hamburg}
\email{schweigert@math.uni-hamburg.de}
\thanks{J.F.\ is supported by VR under project no.\ 621--2003--2385,
and C.S.\ by the DFG project SCHW 1162/1-1.
}

\date{}   

\begin{abstract}
A Morita class of symmetric special Frobenius algebras $A$ in the modular tensor
category of a chiral CFT determines a full CFT on oriented world sheets. 
For unoriented world sheets, $A$ must in addition possess a reversion, i.e.\ an 
isomorphism from \Aop\ to $A$ squaring to the twist. Any two reversions of an
algebra $A$ differ by an element of the group $\Aut(A)$ of algebra automorphisms
of $A$. We establish a group homomorphism from $\Aut(A)$ to the Picard group of 
the bimodule category \CAA, with kernel consisting of the inner automorphisms,
and we refine Morita equivalence to an equivalence relation between algebras 
with reversion. 
         \addforarchive
\end{abstract}

\maketitle

\section{Quantum field theory and categories}

A means for getting to the core of a \qft\ (QFT) is to understand it as a 
functor from some geometric category \calx\ to an algebraic category \cala. 
Since quantum field theory can be analyzed from diverse points of view, 
various such functors, between different types of \cats, have been studied.  

In this paper we consider a specific class of QFT models: two-di\-men\-sio\-nal 
\cfts, or {\em CFTs\/}, for short. Our categorical setup for these can be 
sketched as follows. The objects \X\ of the geometric category \calx\ are 
compact two-di\-men\-sio\-nal manifolds with certain decorations -- 
disjoint labeled (germs of) arcs in the interior of \X\ and/or on the boundary 
$\partial\X$, and/or labeled curves in the interior of \X. The morphisms of 
\calx\ are mapping classes $\varphi{:}\ \X\To\X'$ that are compatible with the 
decorations. The labels for the decorations are taken from data which we 
collectively denote by \CC; accordingly we write $\calx\eq\calxc$. The 
algebraic category \cala\ is the category \Vectp\ of pointed finite-dimensional 
complex vector spaces; its objects are pairs $(W,w)$ consisting of a 
finite-dimensional \complex-vector space and an element $w\iN W$, and its 
morphisms are linear maps. On objects, the CFT functor 
$\cftf{:}\ \calxc\To\Vectp$ acts by mapping a world sheet \X\ to the pair 
  $$
  \cftf(\X) = \big(\, H(\X) \,, \Cf(\X) \,\big) 
  $$
consisting of the space of {\em conformal blocks\/} on \X\ and of the {\em 
correlator\/} of \X.  We will refer to objects \X\ of \calx\ as {\em world 
sheets\/} (a terminology borrowed from string theory) and sometimes slightly 
abuse notation by using the symbol \X\ also for the underlying undecorated 
manifolds.

Deliberately, several basic aspects of this setup, such as the precise form and
physical significance of the decorations and of the data \CC,
have not been specified above. More explanations will be given in due 
time, though for lack of space various details will be suppressed. Before 
doing so it is, however,  wise to examine the following simpler situation. 
Take \calx\ to be the category \calyc\ whose objects are compact oriented closed
two-manifolds with disjoint labeled arcs and whose morphisms are mapping classes 
which map arcs to arcs. The labels for the arcs on objects of \calyc\ are taken 
from some data \C, which are required to form a braided monoidal category: an 
arc is labeled by an object of \C. Further, take \cala\ to be the ca\-te\-go\-ry
\Vect\ of finite-dimensional complex vector spaces. Then we can consider a 
functor
  $$
  \cftco:\quad \calyc \to \Vect \,,
  $$
known as the functor of {\em chiral\/} \cft. (To distinguish the functor \cftf\ 
above from the one of chiral CFT, we call \cftf\ the functor of {\em full} 
\cft.) If it exists, the functor \cftco\ is supposed to be 
determined uniquely by the category \C. When \C\ 
is a {\em \mtc\/} (see section 2), then \cftco\ exists and is indeed 
well-known: it is a two-di\-men\-sio\-nal topological modular functor in the 
sense of \cite{BAki}. To each \mtc\ there is associated such a functor, together
with a three-di\-men\-sio\-nal topological \qft\ \cite{retu2,TUra,BAki}.

A forgetful functor $\Fx{:}\ \calxc\,{\to}\calyc$ is obtained
by viewing world sheets \X\ as pairs 
  $$
  \X = (\Xh,\tau)
  $$
consisting of their {\em double\/} \Xh\ and an orientation-re\-ver\-sing 
involution $\tau$ of \Xh\ (see section 3), and then forgetting $\tau$ and
suitably manipulating the decorations. Together with the
obvious forgetful functor $\Fv{:}\ \Vectp\To\Vect$, $\Fx$ fits into a diagram
  \be
  \begin{array}{ccc}
  \calxc & \stackrel{\scriptstyle\cftf}\verylongrarr  & \Vectp \\[.3em]
  |&&|\\[-1.0em]
  \hspace*{-1.3em}{\scriptstyle\Fx}&&{\scriptstyle\Fv}\hspace*{-1.3em}
  \\[-.6em] \downarrow&&\downarrow \\[-.3em]
  \calyc & \stackrel{\scriptstyle\cftco}\verylongrarr & \Vect
  \eear
  \label{4cf}\ee
of four categories and four functors (all of which are symmetric monoidal).
As we will demonstrate below, the construction of the functor \cftf\ can be
conveniently divided in two steps of which the first amounts to requiring
that this diagram is commutative. 
Thus in particular the spaces of conformal blocks of a full CFT are those 
of the associated chiral CFT, $H((\Xh,\tau)) \eq \cftco(\Xh)$.

We still need to specify the data \CC\ on which \calxc\ depends, as well as the 
corresponding data for the functor \cftf. We must in fact consider two different
categories \calxc, one in which the world sheets are oriented and one in which 
they are not, and accordingly there are two full CFT functors. Here we restrict 
our attention to the oriented case; remarks on the unoriented case will be added
later. For compatibility with \erf{4cf}, the data \CC\ must fit with those of 
\C. Now \C\ is in particular a monoidal category, or what is the same, a 
2-category with a single object. \CC\ generalizes this aspect: it is a 
2-category with precisely two objects. \CC\ has four 1-cells: the monoidal 
category \C, another monoidal category \Cst, and two categories \ML\ and \MR\ 
which are right and left module categories \cite{crfr} over \C, respectively; 
\CC\ must have the further properties that \C\ is braided and that \Cst\ is 
equivalent, as a module \cat\ over \C, to the \cat\ ${\rm Fun}_\C(\ML,\ML)$
of module endofunctors of \ML. Specifying also the functor \cftf\ requires one
additional datum, namely an \alg\ $A$ in \C\ such that \ML\ is equivalent to 
the \cat\ \CA\ of left $A$-modules. Such an \alg\ exists and is determined 
up to Morita equivalence \cite{ostr}, and indeed up to equivalence the whole 
2-category \CC\ can be reconstructed from \C\ and $A$. The algebra $A$ is,
however, only auxiliary; from the correlators obtained with a particular 
choice of $A$, the correlators obtained with any Morita equivalent algebra
$A'$ can be determined uniquely, and indeed the corresponding full CFTs do not
differ in any observable quantity. Nevertheless, for the concrete description 
of the full CFT functor \cftf\ a choice of $A$ within its Morita class must be 
made, and accordingly from now on we denote this functor
by \cfta\ and write $\cfta(\X) \eq (H_\C(\X),\CfA(\X))$.

The rest of this paper is organized as follows. Sections 2 and 3 provide further
information on chiral and full CFT. Some aspects of a `TFT construction' of 
$\CfA(\X)$ are described in sections 4 to 6, while section 7 gives a brief 
outlook to related issues that we do not explain in this paper.
Finally, sections 8 to 11 settle some questions that arise when discussing 
full CFT on unoriented world sheets.

Before proceeding, let us briefly mention other possibilitites for studying
QFT in a categorical framework. In one approach (see e.g.\ \cite{atiy6,sega8}),
which has in particular be discussed for CFTs and for topological \qfts\ (TFTs),
the relevant geometric category \calx\ is a cobordism category, while \cala\ is 
a category of vector spaces. The three-dimensional TFT
variant of this approach, in a formulation closely following the one
of \cite{TUra}, will be used as a tool in our analysis of the functor \cfta.
In another, quite distinct, framework which emphasizes the principle of 
causality \cite{brfv}, the geometric \cat\ has Lorentzian manifolds 
(`space-times') as objects and isometric embeddings as morphisms, while 
\cala\ is the \cat\ of unital C$^*$-\alg s.
Common to all these approaches is that QFT is considered simultaneously on a 
large class of spaces, or space-times.



\section{Chiral conformal field theory and modular tensor categories}

To characterize a specific model of QFT one must in particular have a grasp on
its fields and symmetries. 
The symmetries of a chiral CFT include conformal symmetries. These can be 
encoded in a Virasoro vertex algebra $\V_{\rm Vir}$, which furnishes in 
particular a \rep\ of the Virasoro Lie algebra. But most models have additional 
symmetries, so that the symmetry structure is a larger conformal vertex 
algebra $\V\,{\supseteq}\,\V_{\rm Vir}$ (see e.g.\ \cite{HUan}).
The spaces of fields of the chiral CFT are then given by \V-modules.
The braided monoidal category \C\ that supplies the data for the geometric 
\cat\ of a chiral CFT is the \rep\ \cat\ \repV\ of \V, which has
\V-modules as objects and intertwiners between \V-modules as morphisms. 
The monoidal structure on \repV\ is given by the tensor product $\otimes$ of 
\V-mo\-dules and of intertwiners, the tensor unit $\one$ being \V\ itself 
\cite{hule3.5}. Envoking coherence, we tacitly pass to an equivalent \cat\ 
for which both $\otimes$ and $\one$ are strict.

In the sequel, we will restrict our attention to the case that \V\ is a 
{\em rational\/} conformal vertex algebra in the sense that it satisfies the 
conditions of theorem 5.1 of \cite{huan24}. (That is, \V\ is $C_{2}$-cofinite 
and self-dual as a \V-module and satisfies $\V^{(n)}\eq0$ for $n\,{<}\,0$ 
and $\V^{(0)}\eq\complex\one$, every simple \V-module not isomorphic to \V\,\ 
has positive conformal weight, and every $\mathbb{N}$-gradable weak \V-module 
is fully reducible.) We are then dealing with a (chiral) rational CFT, or 
{\em RCFT\/}, for short. Among the conformal vertex algebras the rational ones 
are distinguished by the fact\,%
 \footnote{~%
 The conditions in the definition of rationality used here can possibly be
 relaxed. Various similar notions of rationality that have been discussed in
 the literature are not sufficiently strong for our purposes, however.}
\cite{hule3.5,huan24} that $\C\eq\repV$ has a number of peculiar properties:
\\[-.8em]
  \myitems
\item[(i)\phantom{ii}]
  The tensor unit is simple.
\item[(ii)\phantom{i}]
  \C\ is abelian, \complex-linear and semisimple.
\item[(iii)]
  \C\ is ribbon:\,%
    \footnote{~%
    Besides the qualifier `ribbon' \cite{retu2}, which emphasizes the
    similarity with the properties of ribbons in a three-manifold, also the
    terms `tortile' \cite{joSt6} and `balanced rigid braided' are in use.}
  There are families $\{c_{U,V}\}$ of braiding, $\{\theta_U\}$ of
  twist, and $\{d_U,b_U\}$ of evaluation and coevaluation morphisms
  satisfying the relevant properties.
\item[(iv)]
  \C\ is Artinian (or `finite'), i.e.\ the number of isomorphism
  classes of simple objects is finite.
\item[(v)\phantom{i}]
  The braiding is maximally non-degenerate: the numerical matrix $s$ with
  entries $s_{i,j} \,{:=}\, (d_{U_j}^{}\oti\tilde d_{U_i}^{}) \cir
  [\id_{U_j^\vee}\oti(c_{U_i,U_j}^{}{\circ}\,c_{U_j,U_i}^{})\oti\id_{U_i^\vee}]
  \cir (\tilde b_{U_j}^{}\oti b_{U_i}^{})\,$ 
  is invertible.
\end{itemize}    

\noindent
Here we denote by $\{U_i\,|\,i\iN\I\}$ a (finite) set of representatives
of isomorphism classes of simple objects; we also take $U_0\,{:=}\,\one$ as the
representative for the class of the tensor unit. A monoidal \cat\ with the 
properties listed above is called a {\em \mtc\/} \cite{TUra}.

It is worth mentioning that every ribbon \cat\ is {\em sovereign\/}, i.e.\
besides the left duality given by $\{d_U,b_U\}$ there is 
also a right duality (with evaluation and coevaluation morphisms to be denoted
by $\{\tilde d_U,\tilde b_U\}$), which coincides with the left duality in the
sense that ${{}^\vee\!}U\eq U^\vee$ and ${{}^\vee\!}\!f\eq f^\vee$.
Below we will make ample use of Joyal\,-\,Street \cite{joSt5} type diagrams
for morphisms of a strict ribbon \cat. For instance, the
ribbon structure morphisms and the entries of $s$ are depicted by
\\
  \inclpic{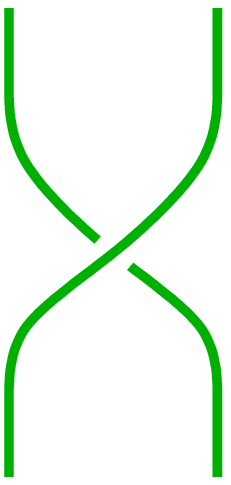}{3} {(105,40)(-36,10) \put(-34,18){$c_{U,V}^{}\,{=}$}
                  \setulen75 \put(4,0){\turlabl U} \put(28,0){\turlabl V} }
  \inclpic{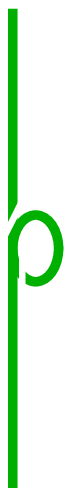}{3} {(56,40)(0,10) \put(-27,18){$\theta_{U}^{}\,{=}$}
                  \setulen75 \put(6,0){\turlabl U} }
  \inclpic{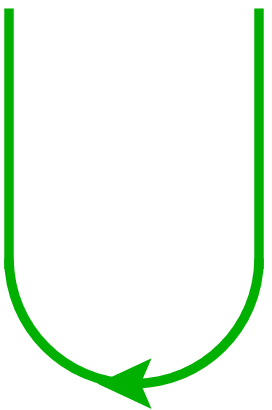}{3} {(72,40)(0,7) \put(-28,15){$b_{U}^{}\,{=}$} \setulen75
                  \put(3,38){\turlabl U} \put(32,38){\turlabl{U^{\!\vee}}} }
  \inclpic{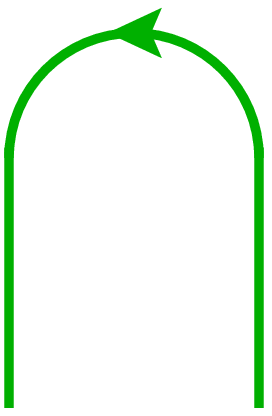}{3} {(73,40)(0,7) \put(-28,15){$d_{U}^{}\,{=}$} \setulen75
                  \put(3,0){\turlabl{U^{\!\vee}}} \put(32,0){\turlabl U} }
  \inclpic{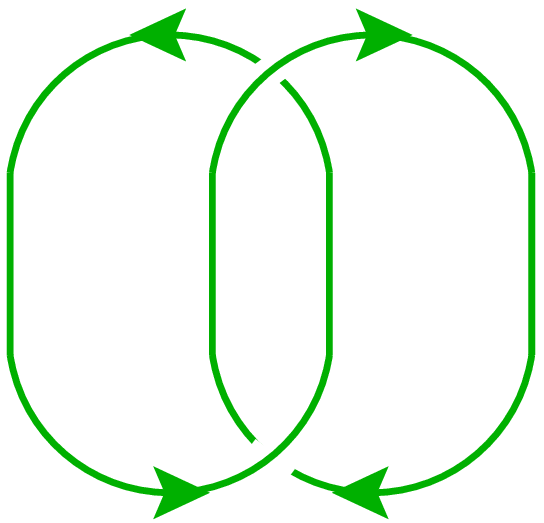}{3} {(44,40)(0,10) \put(-30,18){$s_{i,j}^{}\,{=}$} \setulen50
                  \put(23,53){\Turlabl{U_i}} \put(61,38){\turlabl{U_{\!j}}} } 
\\[.4em]

Let us also remark that other approaches to conformal \qft, not based on vertex
\alg s, exist. In particular one can work with nets of von Neumann \alg s 
instead, see e.g.\ \cite{rehr22}. This setting leads again to \mtcs, albeit with
the extra property that the (quantum) dimensions $\dim(U)\,{:=}\,\tilde d_U
\cir b_U \eq d_U\cir\tilde b_U\iN\complex\,\id_\one$ are real and positive.
We will not be concerned with the issue where the \cat\ \C\ under
study comes from. Moreover, while in sections 3\,--\,7 \C\ will always
stand for a \mtc, various results that we will mention actually remain valid 
in a more general setting. And in sections 8\,--\,11, \C\ will generally not be
assumed to be modular, but essentially only to be an additive \koerper-linear 
ribbon category.


\section{Full conformal field theory}

There exist some physical problems, such as the fractional quantum Hall effect 
(see e.g.\ \cite{fpsw}), for which chiral CFT is relevant. But more often,
e.g.\ in the study of phase transitions, percolation, impurity 
problems, or string theory, it is full rather than 
chiral CFT that matters. Recall that, unlike in chiral CFT, in full CFT the 
objects of \calx\ -- the world sheets \X\ -- can have nonempty boundary. 
In applications, a boundary can arise e.g.\ as a cluster boundary 
in an effectively one-dimensional condensed matter system,
or as the world line of an end point of an open string in string theory.\,%
  \footnote{~%
  Instead of arcs on \X, one may also use parametrized circles to implement
  decorations. Then in addition to the `physical' boundary components
  of \X\ occurring here, one also has `insertion', or `state', boundary 
  components. For details as well as references see \cite{rffs}.}

In physical terminology, an arc on \X\ specifies the location (together with a
germ of local coordinates around it) of a `field 
insertion', while the label of the arc specifies the type of field that
is inserted. As mentioned above, in chiral CFT fields correspond to
\V-modules, and hence arcs are labeled by objects of \C. In contrast, in full 
CFT the label of an arc in the interior of \X\ involves a {\em pair\/} of 
objects of \C, say (taking, without loss of generality the objects to be simple)
$(U_i,U_j)$ with $i,j\iN\I$. Such field insertions are called {\em bulk 
fields\/}. On $\partial\X$ one can have another type of fields, the 
{\em boundary fields\/}, and in the presence of defect lines (labeled curves 
in the interior of \X) there are also {\em defect fields\/}, which
generalize bulk fields. In the sequel we concentrate on the bulk fields.

In terms of the chiral symmetry \V, the prescription to work with a pair of 
objects of \C\ means that bulk fields carry two \rep s of \V; 
or put differently: a \rep\ of two copies of \V. The latter formulation is 
indeed quite suggestive; the two copies are referred to
as ``left- and right-moving'' or ``holomorphic and antiholomorphic'' world 
sheet symmetries, respectively, mimicking the common terminology for the two 
types of solutions to the classical equation of motion of a free boson 
field in two-dimensional Minkowski and Euclidean space, respectively.
(But recall that there are two types of full CFT, with
oriented and unoriented (including in particular unorientable) world sheets,
respectively. The present terminology is appropriate only in the oriented case.)
For oriented full CFT, the total left and right chiral symmetries can actually 
be different. What we denote by \V\ is a sub\alg\ of symmetries that is 
contained both in the left and right chiral symmetries. Since by assumption \V\ 
is a rational conformal vertex algebra, this excludes so-called heterotic RCFTs,
in which rationality is only present when different extended symmetries are
taken into account for the left- and right-moving part.  

When analyzing full RCFT, important tools are supplied by the corresponding
chiral CFT that according to the diagram \erf{4cf} shares the underlying \mtc\ 
\C. In the sequel we take the attitude that this chiral CFT is 
sufficiently well under control, so that the interesting part of discussing 
the full CFT is the particular way it is related to the chiral CFT. 

At the level of geometric categories, the relationship is rather simple: The 
object \Xh\ of \calyc\ to which a world sheet $\X\iN\obj(\calxc)$ gets mapped 
by the forgetful functor \Fx\ is the {\em double\/} \Xh\ of \X.
The world sheet $\X\eq(\Xh,\tau)$ can be obtained from \Xh\ as the quotient by 
an orientation-re\-ver\-sing involution $\tau$. Conversely, the double can be 
recovered from \X\ as the orientation bundle over \X\ modulo identification of 
the two points in the fiber over each point of $\partial\X$ 
\cite{ALgr,bcdcd,fffs3}. One may also regard \X\ as a real scheme; then \Xh\ is 
its complexification, and the involution $\tau$ just implements the action of 
the generator of the Galois group $\mathcal Gal(\complex,\reals)$ \cite{scfu4}.
To give some examples: the double of a closed orientable world sheet \X\ is 
just the disconnected sum $\Xh\eq\X\,{\sqcup}\,{-}\X$ of two copies of \X\ 
endowed with opposite orientation; and both the disk and the real projective 
plane have the two-sphere as their double.

At this point we should mention that world sheets \X\ must also be endowed with 
a conformal structure (and for certain aspects, even with a metric, compare e.g.\ 
\cite{rffs} for details). Analogously, the objects of \calyc\ carry a complex 
structure; the possible choices of complex structure of \Xh\ are restriced by 
the requirement that the involution $\tau$ is anticonformal. For the 
relation between chiral and full CFT studied here, the conformal 
structure on \X\ (and complex structure on \Xh) is inessential.


\section{The connecting three-manifold and topological field theory}

A world sheet \X\ of a full CFT comes with an arc for each bulk field insertion;
each such arc gives rise to two arcs on the double \Xh. 
(In contrast, for boundary field insertions, which are described by arcs on 
$\partial\X$, there is just a single arc on \Xh.) 
To specify the bulk field insertion, the arc on \X\ is labeled by a pair of
objects of \C, say $(U_i,U_j)$ as above. A field insertion
in chiral CFT, on the other hand, requires a single object as label for its arc. 
Thus in order to relate full CFT on $\X\eq(\Xh,\tau)$ to chiral CFT on \Xh\ we 
must label one of the two arcs on \Xh\ by $U_i$ and the other one by $U_j$.
In short, at the level of bulk fields, going from full to chiral CFT affords a 
{\em geometric separation of left- and right-movers\/}.

In the chiral CFT the information that the arc labels $U_i$ and $U_j$ arise from 
one and the same bulk field is ignored. A possibility to retain this information 
is to regard \Xh\ as the boundary of a `fattened' world sheet, which we call the
{\em connecting manifold\/} for \X\ and denote by \MX. \MX\ can be defined as 
the interval bundle over \X\ modulo a certain identification over $\partial\X$,
or equivalently as
  $$
  \MX := \big( \Xh\Times[-1,1] \big) /{\sim}  \qquad{\rm with}\qquad
  ([x,\ort],t)\sim([x,-\ort],-t) \,.
  $$
Then $\partial\MX\eq\Xh$, while \X\ is naturally embedded in \MX\ via
$\io\colon \X\,{\mapsto}\,\X\Times\{t{=}0\}\,{\hookrightarrow}\,\MX$. To relate 
the theories on \X\ and \Xh, it is desirable that also the connecting 
three-ma\-ni\-fold \MX, just like \X\ and \Xh, comes along with some QFT. Since
\MX\ plays only an auxiliary role, that three-di\-men\-si\-onal \qft\ should 
require as little structure on \MX\ as possible -- in physics terminology, it 
should be non-dynamical. 
This is achieved by demanding it to be a {\em topological\/} \qft.

As already mentioned at the end of section 1, a $d$-dimensional topological 
\qft, or {\em TFT\/}, furnishes a specific version of a QFT functor, 
in which the geometric category \calx\ is a cobordism category;\,%
  \footnote{~%
  Owing to a non-trivial action of the mapping class group on the objects $E$
  of the three-dimensional cobordism category, an additional rigidification of
  the category is required. Following \cite{TUra}, here we implement this
  refinement by assuming a choice of Lagrangian subspace of $H_1(E,\reals)$ as
  an auxiliary datum. We refer to \cite{TUra}, as well as e.g.\ to
  \cite{fffs3} or section A.1 of \cite{fjfrs}, for details. An alternative
  treatment is possible in a 2-categorical setting, see \cite{till8,bacR}.}
the target \cat\ \cala\ of a TFT is the category of finite-dimensional
complex vector spaces. Details can be found e.g.\ in \cite{atiy6,TUra,lawr3}.
What is relevant for us is actually a variant, called \C-{\em extended\/} 
three-dimensional TFT; we denote the corresponding functor by \tftc.
The objects of the domain \cat\ of \tftc\ are just those of the domain \cat\ 
\calyc\ of the relevant chiral CFT, i.e.\ compact oriented closed two-manifolds 
$E$ with a chosen Lagrangian subspace of $H_1(E,\reals)$ and
decorated with arcs, which in turn are labeled by objects of \C. The morphisms 
are decorated as well: they are oriented three-manifolds \M\ with embedded 
ribbon graphs, i.e.\ finite collections of disjoint oriented ribbons and 
coupons. Each ribbon either forms an annulus or connects coupons and/or arcs 
on $\partial\M$. The pieces of the ribbon graph are labeled by
data from \C, too: ribbons by objects of \C, and coupons by morphisms of \C,
with domain the tensor product of the objects that label the ribbons entering 
the coupon and codomain the tensor product for the ribbons leaving the coupon. 
(Recall that the tensor product $\otimes$ of \C\ is taken to be strict. 
Also, by strictness of $\one$, ribbons labeled by $\one$ are irrelevant, and 
hence will be regarded as invisible.)

An arc on $\partial\M$ carries the same label as the ribbon beginning or ending
at it. In the situation of our interest, where $\M\eq\MX$, we must in 
addition account for the arcs on $\iox\,{\subset}\,\MX$ in the interior of \M. 
To fit into the TFT picture, such an arc labeled $(U_i,U_j)$ should result 
in a coupon with incoming ribbons labeled by $U_i$ and $U_j$. Let us pretend 
for the moment that no other ribbons are to be attached to such coupons (as 
will be explained in the next section, this is only true for a special class 
of full CFTs). Then we arrive at a TFT description of the bulk fields as
\\
\begin{picture}(210,127)
\put(-1,0){\bPo
    \put(0,0)    {\scalebox{.18}{\includegraphics{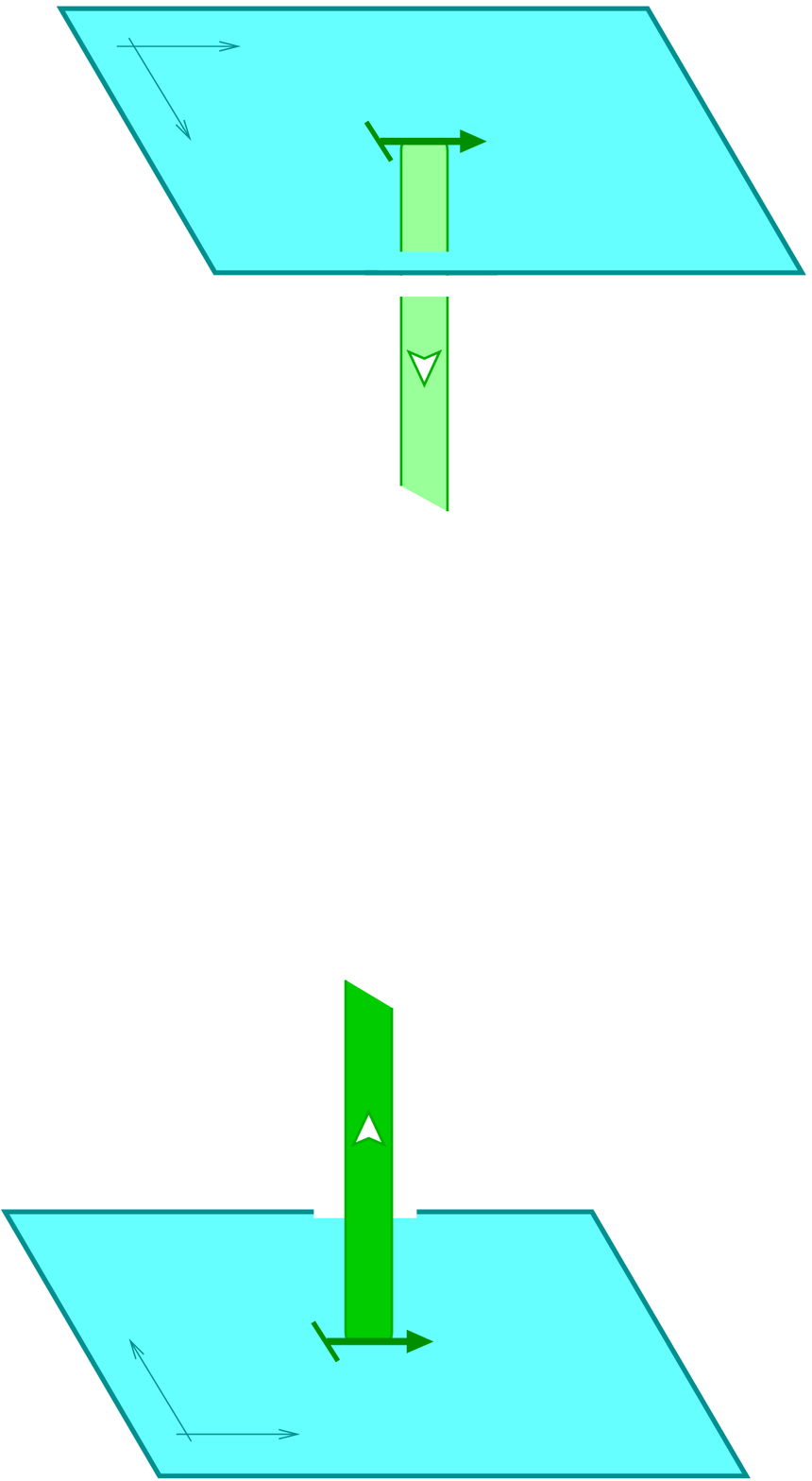}}}
    \put(37,96)  {\small$ U_j $}
    \put(33,7)   {\small$ U_i $}
\eP}
\put(87,0){\bPo
    \put(0,0)    {\scalebox{.18}{\includegraphics{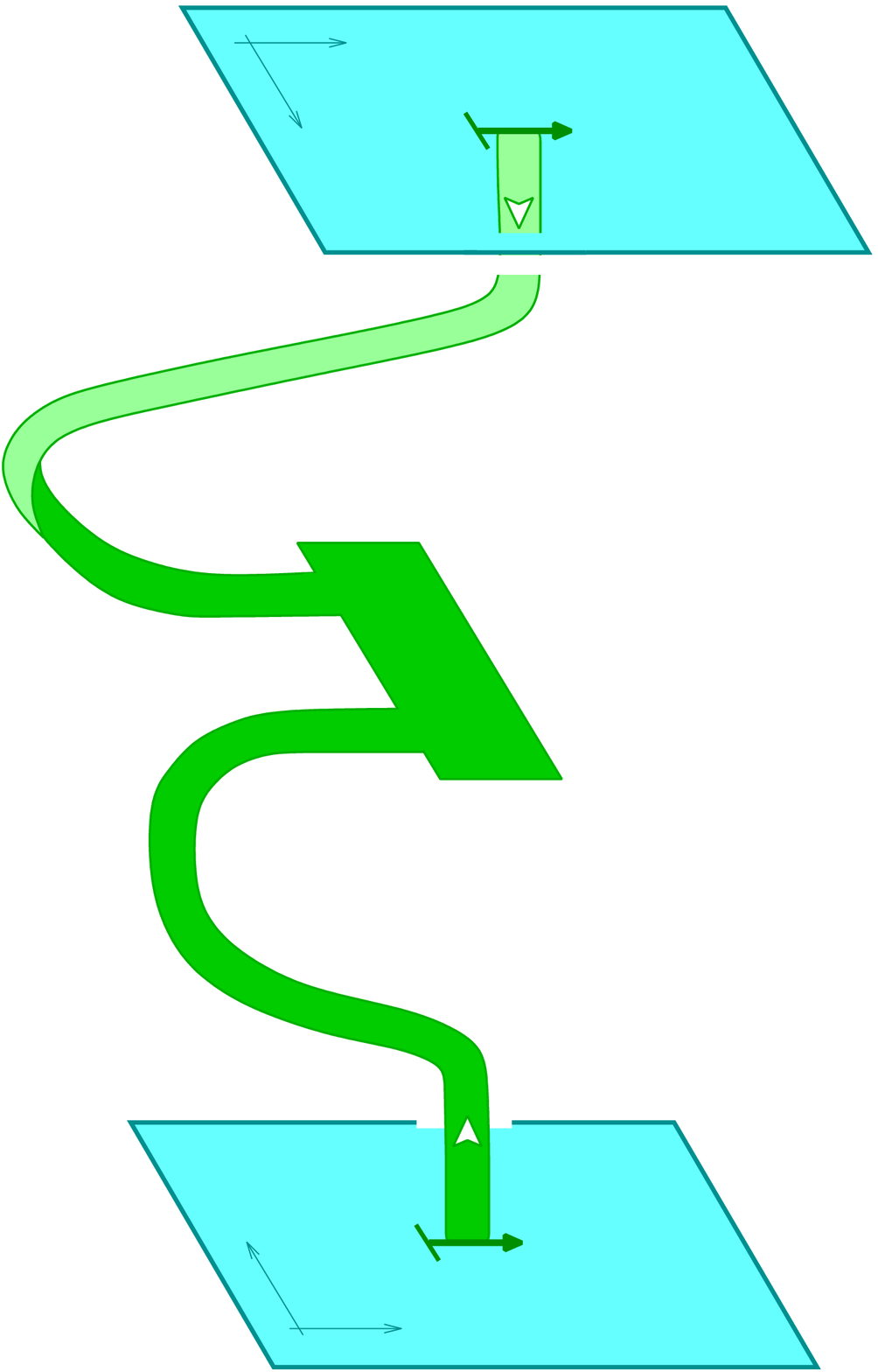}}}
    \setulen82
    \put(0,99)   {\small$ U_j $}
    \put(21.2,45){\small$ U_i $}
    \put(38,67.5){\small$ \alpha $}
\eP}
\put(185,-9){\bPo
    \put(0,0)    {\scalebox{.15}{\includegraphics{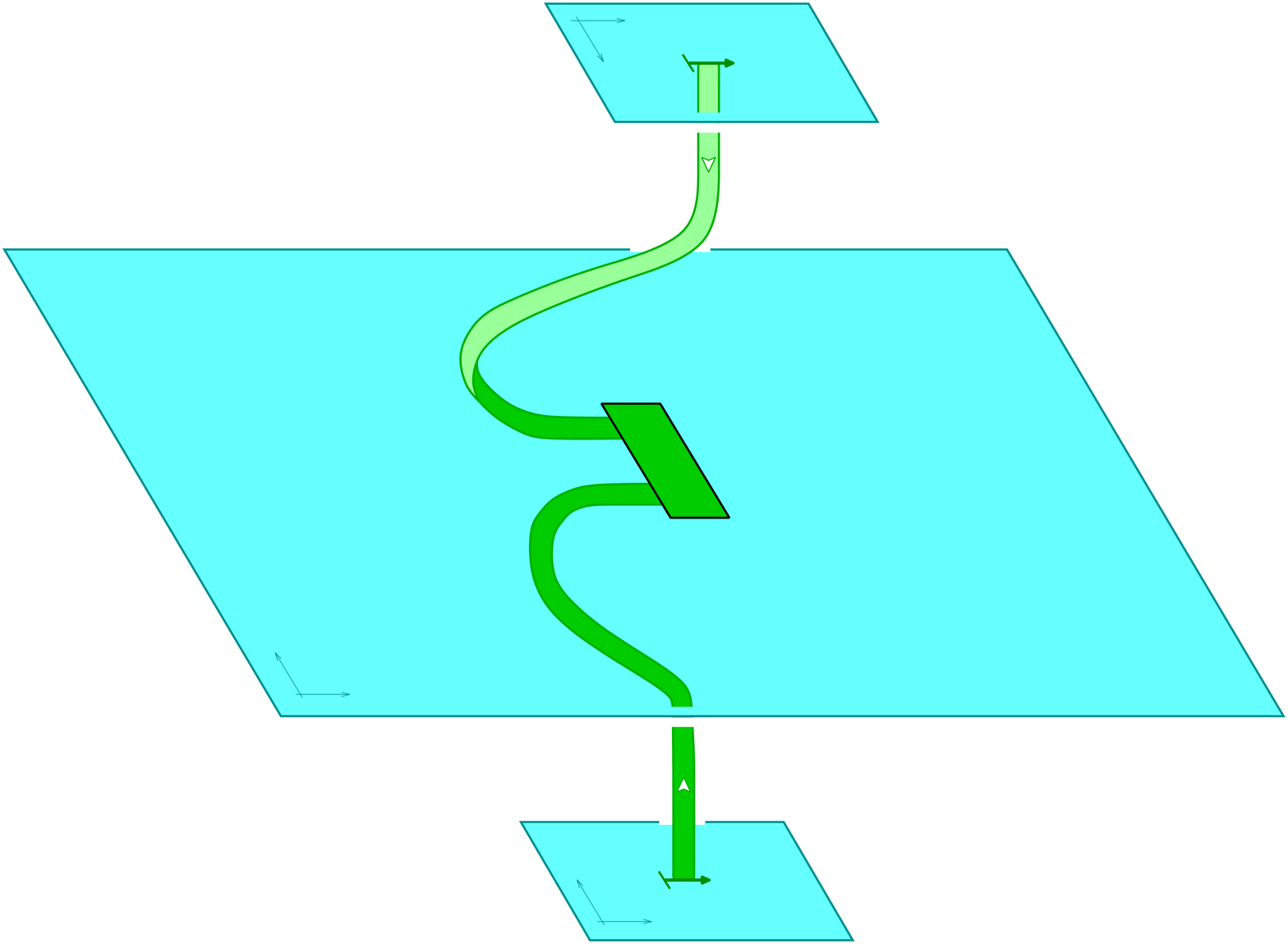}}}
    \setulen60
    \put(133,136) {\small$ U_j $}
    \put(144,69)  {\small$ U_i $}
    \put(146.4,109) {\small$ \alpha $} 
    \put(251,59)  {\small$ \iox $} 
\eP}
\put(25,-15)  {\small (a)}
\put(122,-15) {\small (b)}
\put(247,-15) {\small (c)}
\eP
\\[2.2em]~
The picture (a) just displays the arcs on \Xh\ with emanating ribbons; (b) 
indicates in addition how the ribbons enter the coupon (to be looked at from 
below) in \MX, while (c) shows how this coupon lies in $\iox\,{\subset}\,\MX$.


\section{Full CFT and Frobenius algebras in \C}

As shown in \cite{fffs3}, the description of bulk fields given above can 
indeed be extended to a consistent scheme for constructing a full CFT from its 
underlying chiral CFT with the help of three-dimensional TFT. However, only a 
special class of full CFTs is covered. Indeed, this construction implies in 
particular that the vector space of bulk fields with given chiral labels 
$i,j\iN\I$ is $\Hom(U_i\oti U_j,\one)\,{\cong}\,\delta_{i,j^\vee}^{}\complex$ 
($j^\vee{\in} \,\I$ is the unique label such that $U_{j^\vee}^{}\,{\cong}\,
U_j^\vee$). This is sometimes referred to as $C$-{\em diagonal\/} CFT or, 
when also symmetry preserving boundary conditions are included, as the 
{\em Cardy case\/} of full CFT. In contrast, in general full CFTs these spaces 
can be nonzero for $i\,{\ne}\,j^\vee$, and they can have dimension ${>}\,1$.

In \cite{fuRs4} the missing piece of data was recognized as a certain module
category \ML\ over \C: a full RCFT can be constructed from the corresponding 
chiral CFT (sharing the relevant \mtc\ \C) together with \ML. In short,
  \be
  \mbox{full CFT} ~~=~~ \mbox{chiral CFT}~+~\mbox{module \cat\ \ML\ over \C}
  \,.
  \label{modcat}\ee
A module \cat\ \ML\ over \C\ is a \cat\ \ML\ equipped with a bifunctor
$\ML\Times\C\To\ML$ subject to an appropriate associativity constraint. 
As already mentioned, for any (semisimple) module \cat\ \ML\ there exists 
\cite{ostr} an \alg\ $A$ in \C, unique up to Morita equivalence, such 
that \ML\ is equivalent to the \cat\ \CA\ of left $A$-modules. 

In our application to RCFT, $A$ must be a {\em symmetric special 
Fro\-be\-ni\-us\/} algebra, so that the motto \erf{modcat} may be rephrased 
in a non-Morita invariant manner as
  \be
  \mbox{full CFT} ~~=~~ \mbox{chiral CFT}~+~\mbox{sym.\,sp.\,Frobenius algebra 
  $A$ in \C} \,.
  \ee
A {\em Frobenius\/} \alg\ in \C\ is\,%
  \footnote{~%
  In the \cat\ $\mathcal Vect_\complex$ of finite-dimensional complex
  vector spaces, this definition is equivalent to the more conventional one 
  as an algebra with nondegenerate invariant bilinear form \cite{abra3}.}
a quintuple $A\eq(A,m,\eta,\Delta,\eps)$ such that $(A,m,\eta)$ is a (uni\-tal 
associative) \alg\ and $(A,\Delta,\eps)$ a (counital coassociative) co\alg, 
subject to the compatibility condition that $\Delta$ is a morphism of
$A$-bimodules, i.e. $(\idA\oti m)\cir(\Delta\oti\idA)\eq\Delta
\cir m\eq(m\oti\idA)\cir(\idA\oti\Delta)$; $A$ is called {\em special\/} iff 
$\eps\cir\eta\,{=}$\linebreak[0]$\gamma\,\id_\one$ and $m\cir\Delta\eq\gamma'\,
\idA$ with nonzero complex numbers $\gamma$ and $\gamma'$ 
(implying in particular that $\dim(A)\eq\gamma\gamma'\,{\ne}\,0$); 
and $A$ is called {\em symmetric\/} iff the two isomorphisms 
  \be
  \Phi = ((\eps{\circ}m)\oti\id_{A^\vee}) \cir (\idA\oti b_A) \quad {\rm and}
  \quad \Phi' = (\id_{A^\vee}\oti(\eps{\circ}m)) \cir (\tilde b_A\oti\idA)
  \label{Phi}\ee
in $\Hom(A,A^\vee)$ are equal.

It is worth pointing out that generically \C\ is not symmetric monoidal, and 
as a consequence there are theorems of `braided algebra', such as theorem 
5.20 of \cite{ffrs}, which do not have any substantial classical analogue.

\smallskip

A constructive method for obtaining the full CFT from \C\ and $A$ with the help 
of TFT on the connecting manifold was developed in 
\cite{fuRs4,fuRs8,fuRs10,fjfrs}. It involves in particular a ribbon graph 
$\Gamma$ that is inserted along a triangulation of \iox, with all ribbons (along
edges of the triangulation) labeled by $A$ and all coupons (at the (trivalent) 
vertices of the triangulation of $\X{\setminus}\partial\X$) labeled by the 
product or coproduct of $A$; for details, see e.g.\ appendix A of \cite{fjfrs}.
As for bulk fields, the construction amounts to extending the prescription
of section 4 by introducing additional $A$-ribbons that enter and 
leave the coupon in \iox\ (and are connected to $\Gamma$); schematically,
  \be
\begin{picture}(180,36)(65,31)
\put(18,0){\bPo
    \put(0,0)     {\scalebox{.18}{\includegraphics{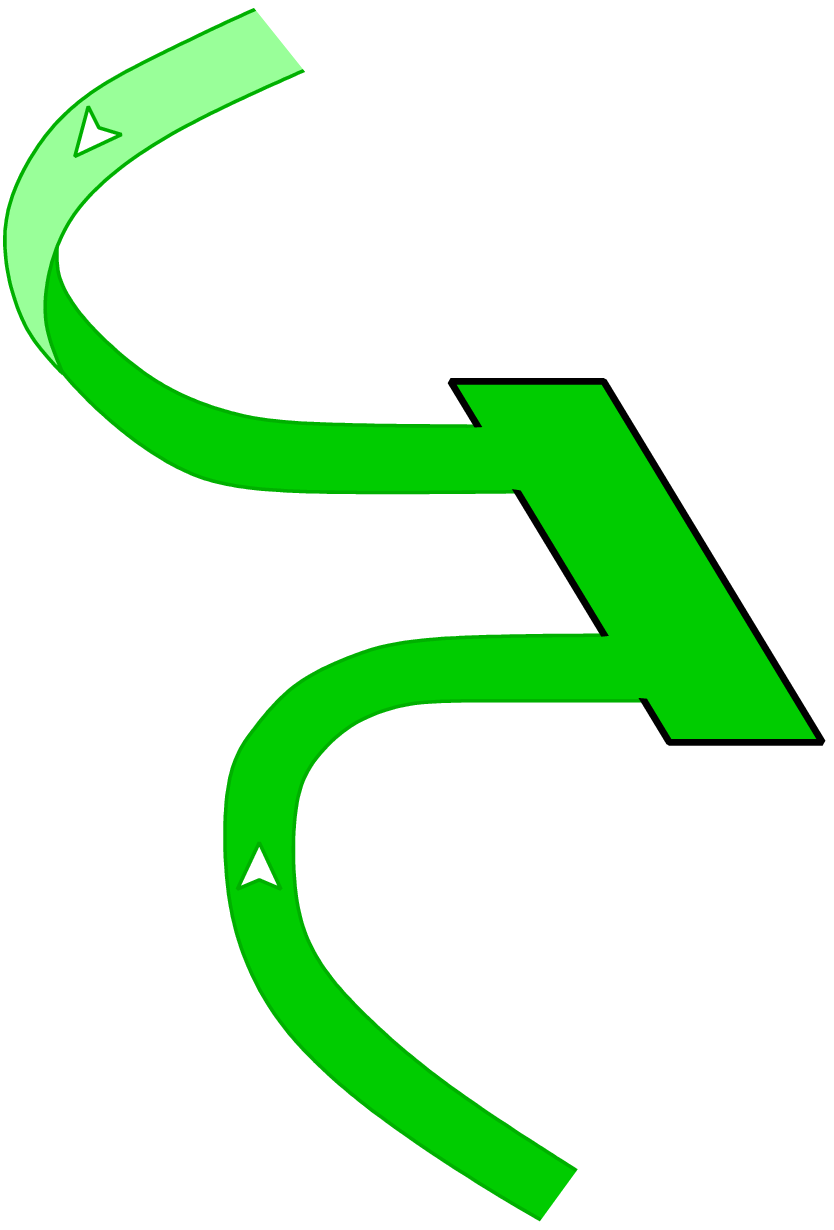}}}
    \put(17,62)   {\small$ U_j $}
    \put(30,35)   {\small$ \alpha $} 
    \put(31,1)    {\small$ U_i $}
\eP}
\put(104,0){\bPo
    \put(0,0)     {\scalebox{.18}{\includegraphics{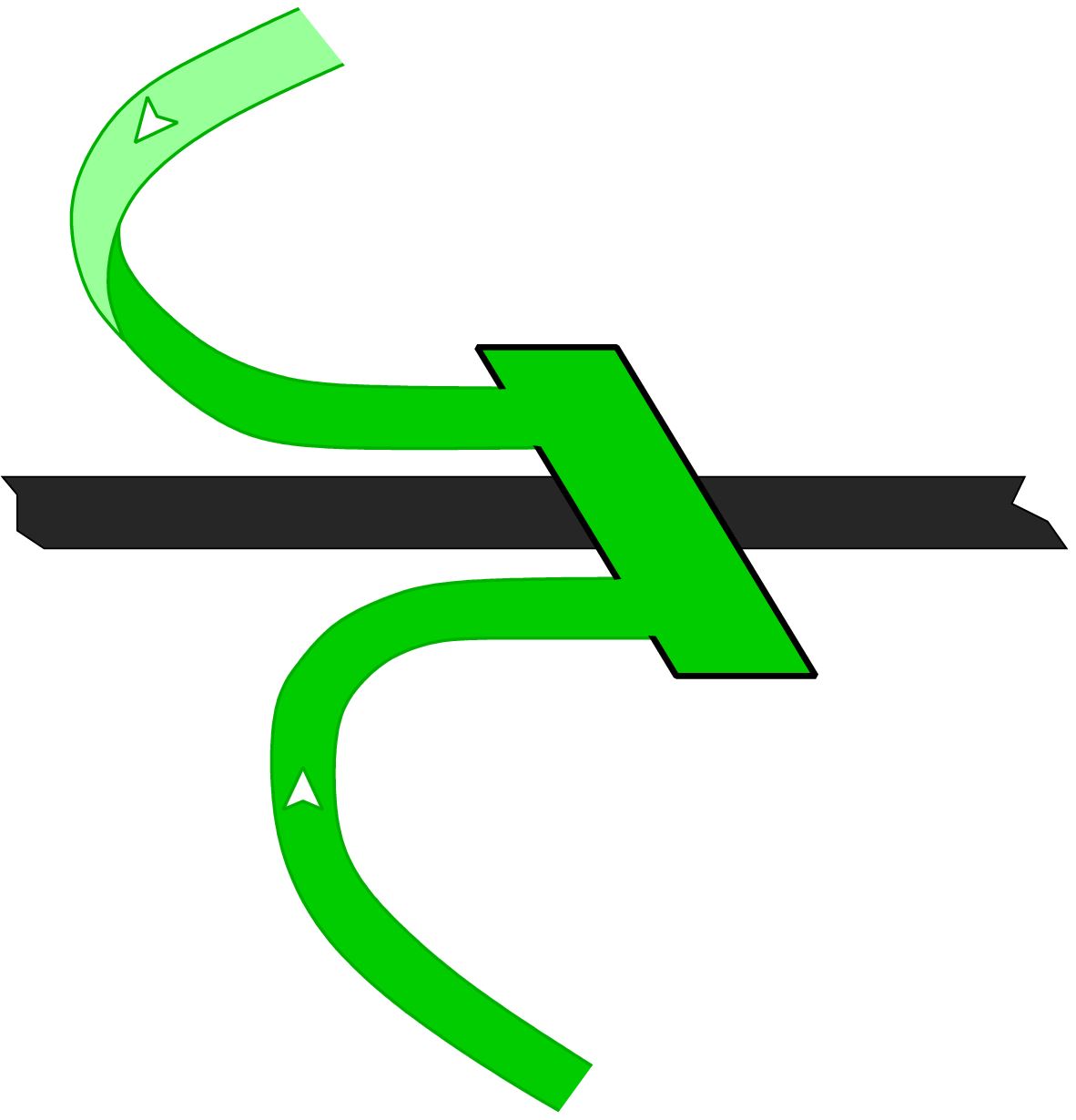}}}
    \put(-7,34)   {\small$ \one $}
    \put(21,62)   {\small$ U_j $}
    \put(34,35)   {\small$ \alpha $} 
    \put(35,1)    {\small$ U_i $}
    \put(63,34)   {\small$ \one $}
\eP}
\put(251,0){\bPo
    \put(0,0)     {\scalebox{.18}{\includegraphics{syd11.eps}}}
    \put(-9,34)   {\small$ A $}
    \put(21,62)   {\small$ U_j $}
    \put(34,35)   {\small$ \alpha $} 
    \put(35,1)    {\small$ U_i $}
    \put(63,34)   {\small$ A $}
\eP}
\put(77,34)   {$\equiv$}
\put(201,33)  {\scalebox{.12}{\includegraphics{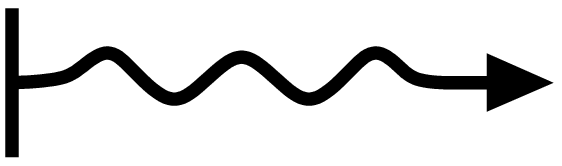}}}
\eP
  \label{bulkrib}
  \ee
~\\[1.9em]
Thus in the general case a bulk field with chiral labels $i,j\iN\I$ is a 
morphism $\alpha\iN\Hom(U_i\oti A\oti U_{\!j},A)$. But to ensure that the 
correlators involving such a bulk field can be nonvanishing, one must in fact 
further restrict to a certain subspace of $A$-bi\-mo\-dule morphisms.
The relevant bimodule structure on the object $U_i\oti A\oti U_{\!j}$
must account for the fact that (referring to picture (c) in the description of
bulk fields in section 4) all $A$-ribbons in \iox\ pass below the $U_i$-ribbon 
and above the $U_j$-rib\-bon, as illustrated in the following picture:
\\
\begin{picture}(80,73)(-161,0)
    \put(-25,0)   {\scalebox{.18}{\includegraphics{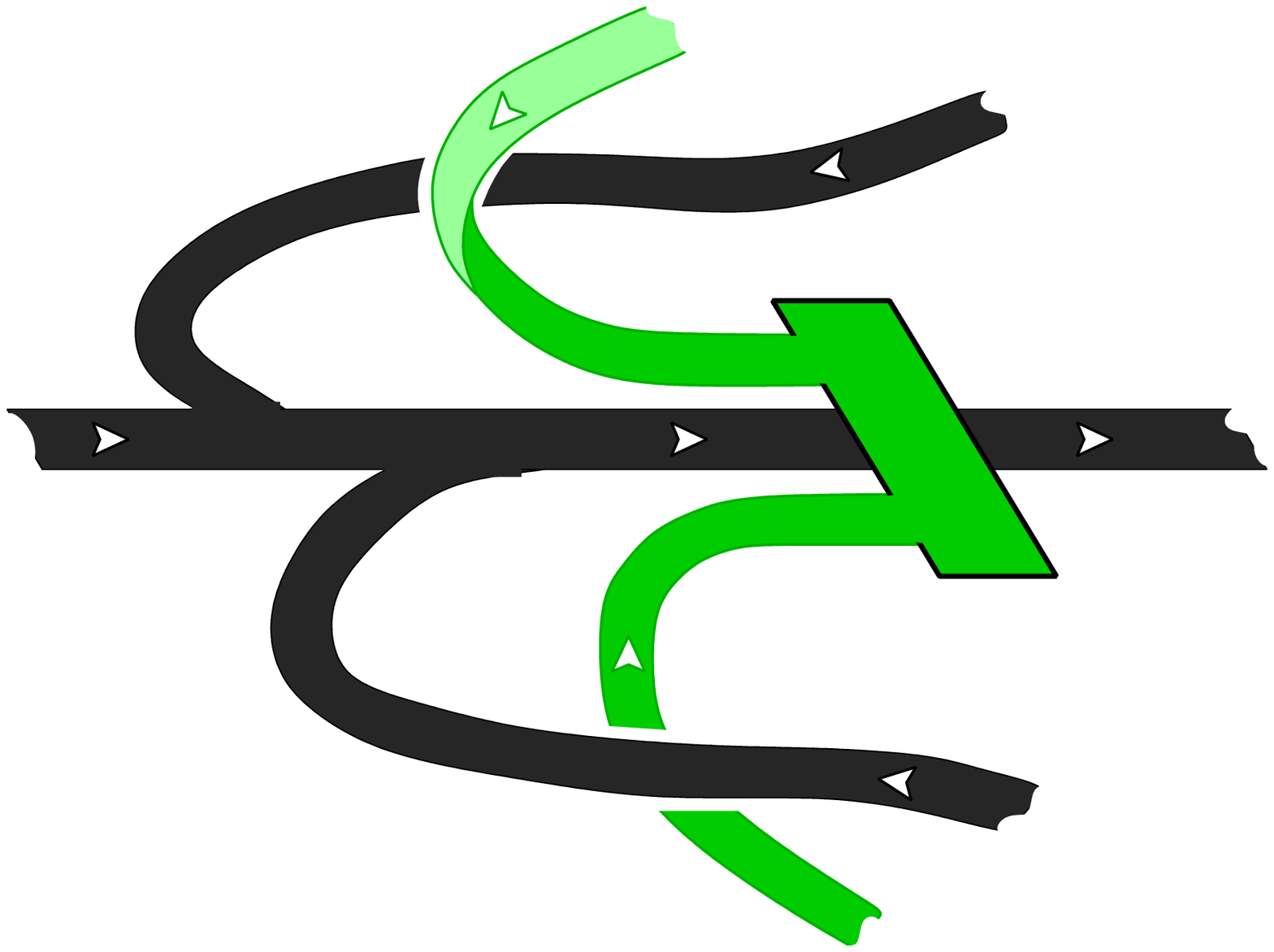}}}
    \put(-32,34)  {\small$ A $}
    \put(24,64)   {\small$ U_j $}
    \put(38,-2)   {\small$ U_i $}
    \put(65,34)   {\small$ A $}
\eP
\\[.6em] 
Accordingly, the commuting
left and right $A$-actions on $U_i\oti A\oti U_{\!j}$ are given by 
$(\id_{U_i}\oti m\oti\id_{U_j}) \cir (c^{-1}_{U_i,A}\oti\idA\oti\id_{U_j})$ 
and $(\id_{U_i}\oti m\oti\id_{U_j})\cir (\id_{U_i}\oti\idA\oti c^{-1}_{A,U_j})$,
respectively. We denote the so defined $A$-bimodule by $U_i\Otip\! A\Otim 
U_{\!j}$. Thus the space of bulk fields with chiral labels $i,j\iN\I$ is 
given by $\HomAA(U_i\Otip\! A\Otim U_{\!j},A)$.

\smallskip

Before proceeding, we mention two pertinent results and some examples of 
\ssfa s. An \alg\ is called {\em simple\/} iff it is simple as a bimodule
over itself, i.e.\ a simple object of the \cat\ \CAA\ of $A$-bimodules.

\begin{proposition}
\label{prop.pertinent}
{\rm(i)}~\,An algebra in \C\ can be endowed with the structure of a \ssFA\ iff
the morphism $\Phi_\natural$, obtained when replacing $\eps$ in formula
{\rm\erf{Phi}} for $\Phi$ by $\xi\epsnat$ with $\xi\iN\complexx$ and 
$\epsnat\,{:=}\,d_A\cir(\idA\oti m)\cir(\tilde b_A\oti\idA)$, is invertible. 
This structure, if it exists, is unique up to the choice of $\xi\iN\complexx$ 
and up to isomorphisms of $A$ as a Frobenius algebra.
\\[.1em]
{\rm(ii)}~Every algebra in the Morita class $[A]$ of a simple \ssFA\ in \C\ is
simple symmetric special Frobenius, too.
\end{proposition}

\begin{proof}
(i) is shown in lemma 3.12 and theorem 3.6(i) of \cite{fuRs4}. 
\\
(ii) follows immediately from proposition 2.13 of \cite{fuRs8}.
\end{proof} 

Concerning examples, first of all the tensor unit $\one$, with all structural 
morphisms identity morphisms, is trivially a \ssFA. Also, for any object $U$ 
of \C, $(U^\vee{\otimes}\,U,\id_{U^\vee}{\otimes}\,\tilde d_U\oti\id_U,
\tilde b_u, \id_{U^\vee}{\otimes}\,b_U\oti\id_U,d_U)$ is a symmetric Frobenius 
algebra, and (since $\eps\cir\eta\eq\dim(U)\,\id_\one$) it is special iff
$\dim(U)\,{\ne}\,0$. These algebras are Morita equivalent to $\one$, 
and hence for each of them the TFT construction gives the $C$-diagonal full CFT. 

A large number of nontrivial examples is provided by the following setup 
\cite{fuRs9}. As a monoidal \cat, the full subcategory \PicaC\ of \C\ whose
simple objects $L$ are the invertible objects of \C\ is determined uniquely
(up to equivalence) by the Picard group \PicC\ of \C\ and a class 
$\psi\iN H^3(\PicC,\complexx)$; $\psi$ specifies the associator of \PicaC.\,%
  \footnote{~%
  Indeed, a categorification of a finite group $K$ is, by definition, a
  monoidal category with Grothendieck ring $\zet K$, and up to equivalence  
  such \cats\ are classified by $H^3(K,\complexx)$. Similarly,
  categorifications of a finite abelian group $T$ are naturally defined as 
  braided monoidal categories, and they are classified \cite{joSt6,FRke} by the
  abelian group cohomology \cite{eima2} $H_{\rm ab}^3(T,\complexx)$.}
Let further $G\,{\le}\,\PicC$ be a subgroup of \PicC\ for which there exists a
class $\omega\iN H^2(G,\complexx)$ such that $\psi_{|G}\eq{\rm d}
\omega$. For each such pair $G,\omega$ there is, up to isomorphism, a unique 
\ssFA\ $A\eq A_{G;\omega}$ in \C, called a Schellekens \alg, for which the 
underlying object is $\bigoplus_{g\in G}L_g$, with $L_g$ a simple object of 
\PicaC\ whose class in \PicC\ is $g$. Moreover, for a given subgroup $G$ the
algebras associated to distinct classes $\omega$ of $H^2(G,\complexx)$ are 
nonisomorphic, and every \ssFA\ $A$ in \C\ all of whose simple subobjects are 
invertible and for which $\Hom(\one,A)\,{\cong}\,\complex$ is of this specific 
form.  

\smallskip

Because of the uniqueness result in part (i) of the proposition, without loss of
generality we can (and do) impose the normalization condition $\eps\cir\eta\eq
\dim(A)\,\id_\one$ (or equivalently, $m\cir\Delta\eq\idA$) on \ssFA s.


\section{TFT construction of RCFT correlators}

One of the major tasks in the study of any \qft\ is to obtain the correlation 
functions, or {\em correlators\/}, of the theory. Roughly speaking, a correlator
is a 
         global
section in a certain vector bundle (involving auxiliary data, 
such as a metric and choices of local coordinates around insertion points) over 
the space of configurations of field insertions. A correlator for the empty 
configuration (no fields at all inserted) is also called a {\em \parfu\/}. In 
full CFT the field configuration for the correlator of $n$ bulk fields is 
specified by the world sheet \X\ and insertion arcs $\gamma_a$ on \X\ with
labels $(U_{i_a},U_{j_a},\alpha_a)$, for $a\iN\{1,2,...\,,n\}$, together with 
further data
which do not concern us here, and similarly for correlators involving boundary
and/or defect fields. These correlators must satisfy three types of consistency 
conditions, known as the (chiral) Ward identities, as modular invariance (or 
locality) and as factorisation constraints, respectively. (For details see e.g.\ 
sections 5.2 and 6.1 of \cite{fuRs10} and sections 4.2 and 5 of \cite{fuRs7},
as well as \cite{rffs}.)

In the construction of \cite{fuRs4,fuRs8,fuRs10,fjfrs}, for brevity to be 
referred to as the {\em TFT construction\/}, any RCFT correlator 
is obtained as the invariant of a three-manifold with embedded ribbon graph
that is built from the data just mentioned. In more detail, the correlator 
$\CfA(\X)$ of the full CFT with \mtc\ \C\ and \ssFA\ $A$ for a world sheet 
$\X\eq(\Xh,\tau)$ 
(with field insertions, which we suppress in the notation) is the vector
  \be
  \CfA(\X) := \tftc(\M_{\emptyset,\Xh})\,1 ~\in \tftc(\Xh)
  \label{CX}
  \ee
in the TFT state space of the double \Xh\ of the world sheet.
Here $1\iN\tftc(\emptyset)\eq\complex$, and $\M_{\emptyset,\Xh} \eq \MX$
is the connecting manifold, regarded as a cobordism from the empty set to \Xh,
with an embedded ribbon graph. This ribbon graph, in turn, is obtained by
assembling various building blocks, such as a fragment of the form
displayed in the last graph of \erf{bulkrib} for each bulk field insertion
and the $A$-colored ribbon graph $\Gamma$ along a triangulation of \iox.
For the entire prescription, see appendix A of \cite{fjfrs}.
That $\CfA(\X)$ obeys the Ward identities is already guaranteed by the fact
that it is an element of $\tftc(\Xh)$; for the particular vector \erf{CX} in 
$\tftc(\Xh)$, the modular invariance and factorisation constraints are 
satisfied as well.

When \X\ is oriented and its boundary $\partial\X$ is empty, so that
$\Xh\eq\X\,{\sqcup}\,{-}\X$, one has
$\tftc(\Xh)\eq\tftc(\X)\,{\otimes_\complex}\,\tftc{(\X)}^*$. According to 
\erf{CX} the correlator is then an element in the tensor product of a space 
for `left-movers' and one for `right-movers'.  This property is known as the
{\em holomorphic factorisation\/} of correlators \cite{witt39}.

The description of correlators as 
     global sections
requires that $\CfA(\X)$ is a single-valued function of the field insertion 
points and of the moduli of \X. That we specify here $\CfA(\X)$ instead as 
an element of some vector space complies with this requirement upon
realizing that vector space as a space of (multivalued) functions, the
space of {\em conformal blocks\/} (also known as chiral blocks).\,%
  \footnote{~%
  In {\em chiral\/} CFT one actually does not have correlators in the sense 
  used here. Rather, there are only the vector spaces of conformal blocks,
  which constitute the fibers of (generically) nontrivial bundles over the
  relevant space of field configurations.}
That the TFT state spaces can indeed be identified with the spaces of 
conformal blocks of the corresponding CFT is still an assumption for general 
RCFT models, but has been established for several important classes such
as the Wess\hy Zumino\hy Witten models.


\section{Dictionary}

Among the results established through the TFT construction,
we mention here only that factorisation and modular invariance
properties of the construction were proven in \cite{fjfrs};
that various properties of the most interesting partition functions were
established in section 5 of \cite{fuRs4} and section 3 of \cite{fuRs8};
that the case of Schellekens algebras is treated in \cite{fuRs9};
and that when expressing $\CfA(\X)$ for a few particular world sheets \X\
(the `fundamental correlators') in certain standard bases of conformal blocks,
one obtains the structure constants of various `operator product expansions',
see section 4 of \cite{fuRs10}.

Another important fact is that, in accordance with \erf{modcat}, for oriented 
\X\ the assignment of a suitably normalized correlator $\CfA(\X)$ to \X\
actually depends only on the Morita class of the \ssFA\ $A$; this has been
discussed e.g.\ in section 4.1 of \cite{fuRs4} and section 1.4 of \cite{fuRs8},
and will be proven in \cite{ffrs5}. Also, as already pointed out in section 1, 
a Morita invariant formulation leads naturally to a setup in terms of a
2-category \CC\ that has precisely two objects. The 1-cells of \CC\ are 
\C, $\ML\,{\cong}\,\CA$, $\MR\,{\cong}\,\CAr$ and $\Cst\eq{\rm Fun}_\C(\ML,\ML)
\,{\cong}\,\CAA$; see \cite{muge8} and section 4 of \cite{ffrs4}.

Instead of providing any further details, we content ourselves to illustrate
some of these aspects in a brief dictionary below.

\begin{center}
\begin{tabular}{lcl}
\phantom.\\[-.5em]\hline\\[-.8em]
\multicolumn3c{\sc Dictionary} \\
\multicolumn3c{between physical concepts and mathematical structures}
\\[.4em]
\hline\\[-.7em]
\dictentry {chiral label}              {object $U\iN\obj(\C)$}
\dictentre {full CFT on oriented \X}   {Morita class $[A]$ of}
\dictentry {}                          {\hspace{1.3em}\ssFA s in \C}
\dictentry {full CFT on unoriented \X%
           \hspace{-.8em}}             {Jandl algebra $(A,\sigma)$ in \C}
\dictentre {space $\{\vPhi_{ij}\}$ of bulk fields}  {vector space
                                        $\HomAA(U_i\Otip\!A\Otim U_{\!j},A)$}
\dictentry {}                          {\hspace{7.7em}of bimodule morphisms}
\dictentry {boundary condition}        {$A$-module $M \iN\obj(\CA)$}
\dictentre {space $\{\vPsi_i^{MM'}\}$} {vector space $\HomA(M\oti U_i,M')$}
\dictentry {\hspace{3.12em}of boundary fields}{\hspace{7.7em}of module morphisms}
\dictentry {defect line}               {$A$-$B$-bimodule $Y\iN\obj(\CAB)$}
\dictentre {space $\{\vTheta_{ij}^{YY'}\}$}  {vector space
                                        $\HomAB(U_i\Otip Y\Otim U_{\!j},Y')$}
\dictentry {\hspace{3.12em}of defect fields}{\hspace{7.7em}of bimodule morphisms}
\dictentry {simple current model}      {Schellekens algebra $A\,{\cong}\,
                                      \bigoplus_{g\in G}L_g$, $G\,{\le}\,\PicoC$}
\dictentry {internal symmetries}       {Picard group\, \PicCAA}
\dictentre {Kramers-Wannier like}      {duality bimodules $Y\iN\obj(\CAA)$\,: }
\dictentry {\hspace{5.9em}dualities}   {\hspace{3.4em}bimodules obeying
                                        $Y^{\vee}\Ota Y\iN\PICCAA$}
\phantom.\\[-.8em]\hline\phantom.
\end{tabular}
\end{center}

 \noindent
For the notions that occur in this list without having been introduced before,
we mainly refer to the cited literature, and only explain the one 
of a Jandl algebra.
\begin{definition}
\label{jandl}
A {\em Jandl algebra\/} $(A,\sigma)$ in \C\ is a \ssFA\ $A$ in \C\ together 
with a morphism $\sigma \iN \Hom(A,A)$ satisfying
  \be
  \sigma \circ \eta = \eta \, , \qquad
  \sigma \circ m = m \circ c_{A,A}^{} \circ (\sigma\oti\sigma)
  \, , \qquad
  \sigma \circ \sigma = \theta_{\!A} \,.
  \label{sig}\ee
\end{definition}
In words, $\sigma$ is an \alg\ isomorphism from the opposite algebra \Aop\ to 
$A$ that squares to the twist; we call it a {\em reversion\/} on $A$.


\section{Automorphisms of Frobenius algebras}

There is obviously a flaw in the above dictionary: Just like for full
CFT on ori\-en\-ted world sheets a Morita class rather than a single
algebra is relevant, also in the non-oriented case the relevant
datum should be a suitable class of Jandl algebras, not an individual one.
The rest of this paper is devoted to structures which are needed to
remedy this flaw. Since some of them are of interest in their own right,
independently of the application to CFT, from this point on \C\ will no
longer be required to be a \mtc. Rather, \C\ is only assumed to be a small 
additive idempotent complete strict ribbon category enriched over 
$\mathcal Vect_\koerper$, with \koerper\ a field, and with the tensor unit 
$\one$ being simple and absolutely simple, i.e.\ 
$\End(\one)\eq\koerper\,\id_\one$. $A$ will again be a \ssfa\ in \C.  

We first study \auto s of $A$.

\begin{definition}
A (unital) {\em algebra \auto\/} of an algebra $(A,m,\eta)$ is an isomorphism
$\varphi\iN\End(A)$ such that $m \cir (\varphi\oti\varphi) \eq \varphi\cir m$ 
and $\varphi\cir\eta \eq \eta$.
\end{definition}

By composition, the algebra \auto s of $A$ form a group, denoted by $\Aut(A)$.
One may analogously define a (co-unital) co-algebra \auto\ of a co-algebra 
$(A,\Delta,\eps)$ as an isomorphism $\varphi\iN\End(A)$ satisfying
$(\varphi\oti\varphi)\cir\Delta\eq\Delta\cir\varphi$ and
$\eps\circ\varphi\eq\eps$, and a Frobenius \auto\ of a Frobenius \alg\ as
one that is both an algebra and a co-algebra \auto.
But when $A$ is symmetric special Frobenius, there is no need to distinguish
between these concepts:

\begin{lemma}
Let $A$ be a \ssfa. A morphism in $\End(A)$ is an algebra \auto\ of $A$ iff it 
is a co-algebra \auto\ of $A$.
\end{lemma}  

\begin{proof}
We show that any algebra \auto\ $\varphi$ of a \ssfa\ $A$ is also a
co-algebra \auto. The opposite implication follows by analogous arguments,
in which the role of the \alg\ and co-\alg\ structures are interchanged.
\\
First, by the defining property of $\varphi$ and sovereignty of \C\ it follows 
that the morphism
$\epsnat\eq d_A\cir(\idA\oti m)\cir(\tilde b_A\oti\idA)\iN\Hom(A,\one)$
obeys $\epsnat\cir\varphi\eq\epsnat$. Invoking proposition \ref{prop.pertinent}
we can set $\eps\eq\xi\epsnat$; hence $\eps\cir\varphi\eq\eps$.
Further, using this result together with the defining property of $\varphi$,
one shows that the isomorphism $\Phi\iN\Hom(A,A^\vee)$ introduced in \erf{Phi}
satisfies $\Phi\cir\varphi\eq\big(\varphi^{-1}\big)^{\!\vee}_{}{\circ}\,\Phi$.
Now the coproduct of $A$ can be expressed through $m$ and $\Phi^{-1}$,
see formula (3.40) of \cite{fuRs4}; when combined with that expression, the 
equality just obtained is easily seen to imply that 
$(\varphi\oti\varphi)\cir\Delta\eq\Delta\cir\varphi$. 
\end{proof}

Next we introduce inner \auto s; to this end we need the \koerper-vector space 
  $$
  \Atop := \Hom(\one,A) \,.
  $$
\Atop\ is a \koerper-\alg\ with unit element $\eta$ and multiplication $\alpha
\,{\multo}\,\beta\,{:=}\,m\cir(\alpha\oti\beta)$ for $\alpha,\beta\iN\Atop$.
If $\alpha$ is invertible with respect to this product, we write $\alpha^{-1}$ 
for its inverse, so that 
$\alpha\,{\multo}\,\alpha^{-1} \eq \eta \eq \alpha^{-1}{\multo}\,\alpha$,
and we denote by \Atopx\ the group of invertible elements of \Atop.
For $\alpha\iN\Atopx$ we set
  \be
  \inn\alpha := m \cir (m \oti \alpha^{-1}) \cir (\alpha\oti\idA) \,.
  \label{inn}\ee

We also introduce the endomorphims
\\
  \begin{picture}(0,74)(-110,0)
    \put(0,0)    {\scalebox{.19}{\includegraphics{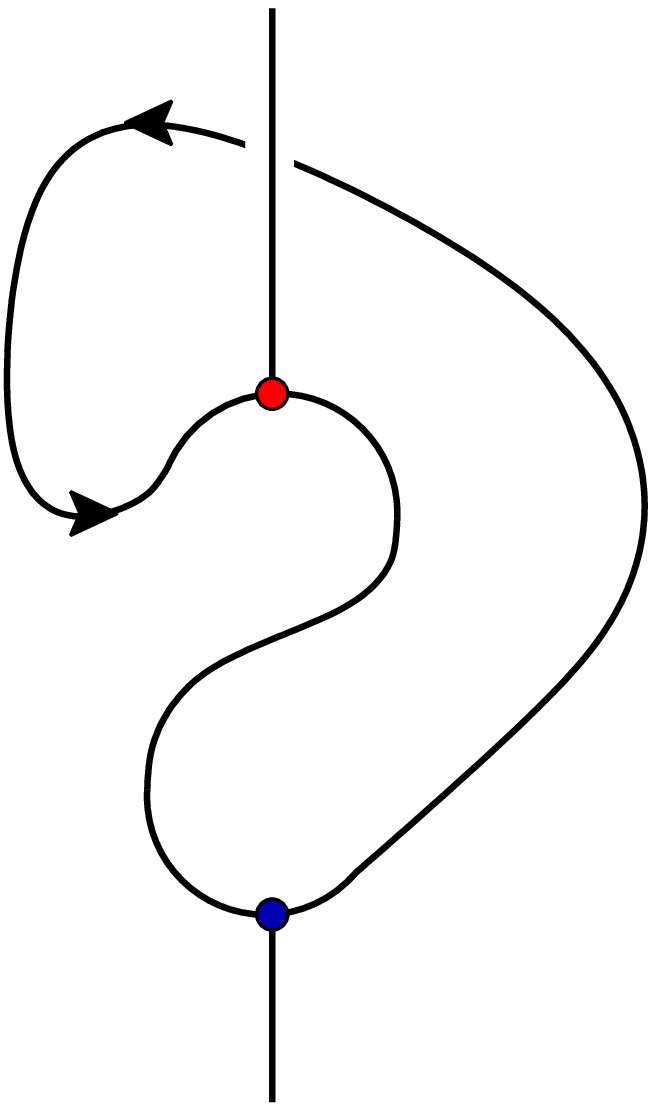}}}
    \put(140,0)  {\scalebox{.19}{\includegraphics{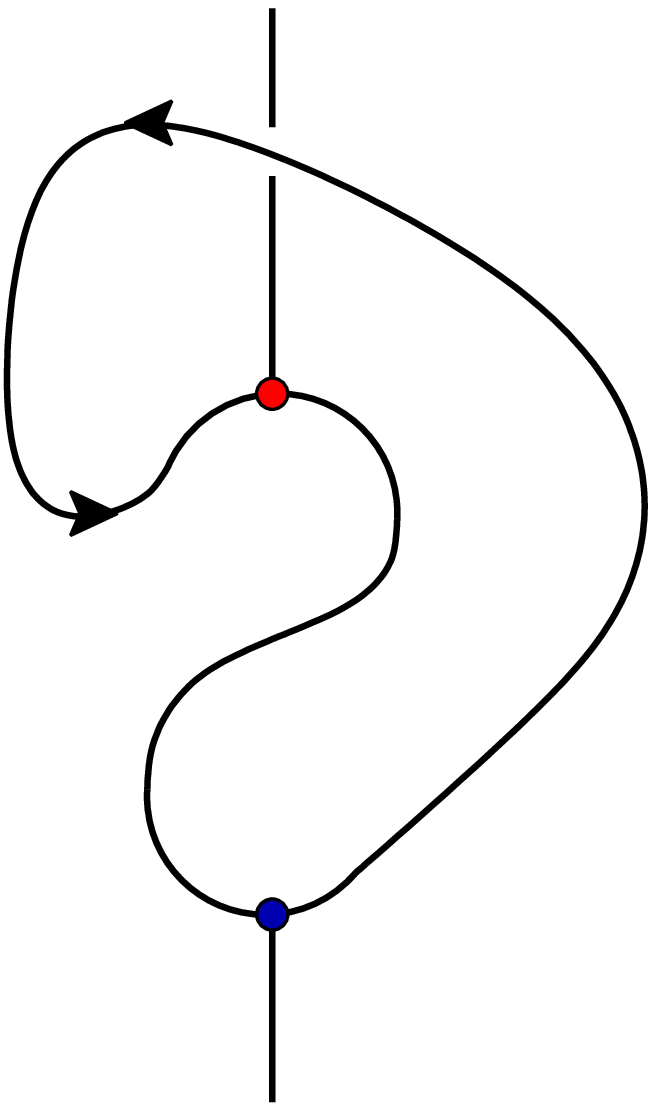}}}
  \setulen50
    \put(-62,55.5) {$P_{\!l}^{}\;:=$}
    \put(1,30)     {\small $A$}
    \put(21.5,-17) {\small $A$}
    \put(23,125.7) {\small $A$}
    \put(73.5,61)  {\small $A$}
    \put(126,55.5) {and}
  \put(279,0){\bPo
    \put(-62,55.5) {$P_{\!r}\;:=$}
    \put(21.5,-17) {\small $A$}
    \put(23,125.7) {\small $A$}
    \put(73.5,61)  {\small $A$}
    \put(1,30)     {\small $A$} \eP}
  \eP
\\[.2em]
of $A$, where 
  \begin{picture}(20,25)(-3,6)
  \put(0,0)    {\scalebox{.19}{\includegraphics{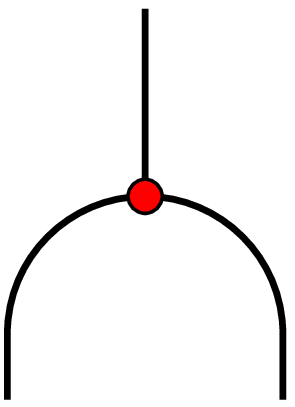}}} \eP
$\iN\Hom(A\Oti A,A)$ is the product and
  \begin{picture}(20,25)(-3,8)
  \put(0,0)    {\scalebox{.19}{\includegraphics{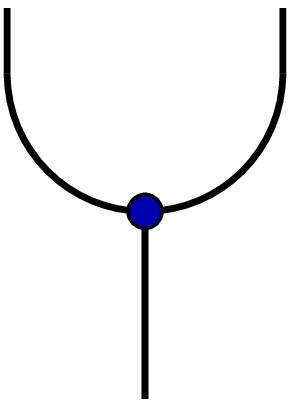}}} \eP
$\iN\Hom(A,A\Oti A)$ is the coproduct. 

One checks that $P_{\!l}^{}\cir\gamma\eq P_{\!r}\cir\gamma$ for 
$\gamma\iN\Atop$. (Also, $P_{\!l}^{}$ and $P_{\!r}$ are idempotents; their
images $C_l(A)\,{:=}\,\Im P_{\!l}^{}$ and $C_r(A)\,{:=}\,\Im P_{\!r}$ are
the left and right {\em centers\/} of $A$, see definition 2.31 of \cite{ffrs}.
A few further properties of $P_{\!l/r}^{}$ are given in lemma 2.29 of 
\cite{ffrs}.) We write
  \be
  \Ctop := \{ \gamma\iN\Atop \,|\, P_{\!l}^{}\cir\gamma\eq\gamma \}
  \ee
(by lemma 2.29(iii) of \cite{ffrs}, \Ctop\ is a \koerper-subalgebra of \Atop)
and denote by \Ctopx\ the group of invertible elements of \Ctop.

\begin{lemma}
\label{lem.div.inn}
The mapping $\alpha\,{\mapsto}\,\inn\alpha$
is a group homomorphism from \Atopx\ to $\Aut(A)$, with kernel \Ctopx.
Thus in particular \Ctopx\ is a normal subgroup of \Atopx.
\end{lemma}

\begin{proof}
By associativity of $A$ one has
$(\alpha{\multo}\beta)^{-1}\eq\beta^{-1}{\multo}\,\alpha^{-1}$, which
implies that $\inn\alpha \cir \inn\beta \eq \inn{\alpha\multo\beta}$.
In particular, $\inn\eta\eq\idA$, and $\inn{\alpha^{-1}}$ is inverse to 
$\inn\alpha$. That for any $\alpha\iN\Atopx$ the morphim $\inn\alpha$ is 
an algebra \auto\ of $A$ then follows by just using associativity of $A$ 
and the unit property of $\eta$.
Finally, again by associativity and the unit property, the equality 
$\inn\gamma\eq\idA$ implies that 
$m\cir(\idA\oti\gamma) \eq m\cir(\gamma\oti\idA)$, which in turn is
equivalent to $\gamma\iN\Ctopx$ (recall that $C_l(A)\,{:=}\,\Im P_{\!l}^{}$\,).
\end{proof}

A slight extension of the consideration in the last part of the proof also
shows that \Ctop\ is contained in the center of \Atop.
Furthermore, by definition of \Ctop\ and of the left and right centers 
$C_{l/r}(A)$ one has $\dimk\,\Ctop\eq\dimk\,\Hom(\one,C_{l/r}(A))$.

We call the algebra \auto s of $A$ that are of the form \erf{inn} the
{\em inner \auto s\/} of $A$, and denote by
  \be
  \Inn(A) := \{ \inn\alpha \,|\, \alpha\iN\Atopx \}
  \ee
the group of inner \auto s. Note that, by lemma \ref{lem.div.inn},
  $$
  \Inn(A) \cong \Atopx / \Ctopx .
  $$


\section{Automorphisms and the Picard group of \,\CAA}

The inner \auto s play a special role when one twists the action of $A$ on 
itself. Owing to associativity and the defining property of an \alg\ \auto, 
twisting the action of $A$ on a left- or right-module by elements 
$\varphi,\psi\iN\Aut(A)$ yields
another left- or right-module structure on the same object. In particular,
  \be
  \aAa\varphi\psi
  := (\,A \,, m\cir(\varphi\Oti\idA) \,, m\cir(\idA\Oti\psi)\,) 
  \label{aAa}\ee
defines an $A$-bimodule structure on the object $A$. Not all of these
bimodule structures are new, though (compare \cite[p.\,112]{vazh}):

\begin{lemma}
\label{lem.aAa}
{\rm(i)}~~For all $\varphi,\,\psi,\,\gamma\iN\Aut(A)$ we have
  $$
  \quad\,\ \aAa\varphi\psi \,\cong\, \aAa{\gamma\varphi}{\gamma\psi}
  \quad \mbox{\em as $A$-bimodules.}
  $$
{\rm(ii)}~\,For all $\varphi,\,\psi,\,\varphi',\,\psi'\iN\Aut(A)$ we have,
as $A$-bimodules,
  $$
  \aAa\varphi\psi \otA \aAa{\varphi'\!}{\psi'}
  \,\cong\, \aAa{\psi^{-1\!}\varphi}{\varphi'{}^{-1\!}\psi'} \,.
  $$
~\\[-.9em]
{\rm(iii)}~For any $\varphi,\,\psi\iN\Aut(A)$,
$\aAa\varphi\psi$ is an invertible object of \CAA.
\\[.4em]
{\rm(iv)}~We have\hspace{5.2em}
  $
  \aAa\idss\psi \cong \aAa\idss\idss \;\Longleftrightarrow\; \psi\iN\Inn(A) \,.
  $
\end{lemma}

\begin{proof}
(i)~~Using that $\gamma$ is an \alg\ \auto, one easily sees that $\gamma$ 
intertwines both the left and the right action of $A$ on $\aAa\varphi\psi$ and 
on $\aAa{\gamma\varphi}{\gamma\psi}$.
Since $\gamma$ is invertible, this establishes the claim.
\\[.3em]
(ii)~\,The statement follows from the fact that
$\aAa\varphi\psi \otA \aAa{\varphi'\!}{\psi'}\,{\cong}\,
\aAa{\psi^{-1\!}\varphi}\idss \otA \aAa\idss{\varphi'{}^{-1\!}\psi'}$,
which in turn is a direct consequence of (i).
\\[.3em]
(iii)~By (i) and (ii) we have in particular
$\aAa\varphi\psi \otA \aAa\psi\varphi \,{\cong}\, \aAa{\psi{}^{-1\!}\varphi}
{\psi{}^{-1\!}\varphi} \,{\cong}\, \aAa\idss\idss \eq A$. 
\\[.3em]
(iv)~Let $\psi\eq\inn{\alpha^{-1}}\iN\Inn(A)$. Then 
$m\cir(\idA\oti\alpha)\iN\End(A)$ (which is invertible, the inverse being
$m\cir(\idA\oti\alpha^{-1})$) intertwines the right $A$-action on the bimodules
$\aAa\idss\psi$ and $\aAa\idss\idss$ (and, trivially, the left $A$-action as
well). Thus $\aAa\idss\psi$ and $\aAa\idss\idss$ are isomorphic as 
$A$-bimodules.
\\
Conversely, let $\varphi\iN\End(A)$ be an intertwiner between
$\aAa\idss\psi$ and $\aAa\idss\idss$, i.e.\ satisfy
\\
\begin{picture}(20,89)(-34,0)
\setulen79
  \put(0,0)       {\scalebox{.36}{\includegraphics{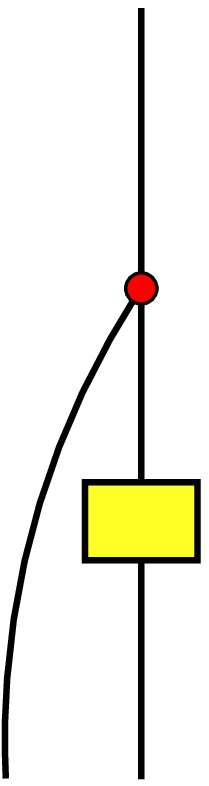}}}
  \put(84,0)      {\scalebox{.36}{\includegraphics{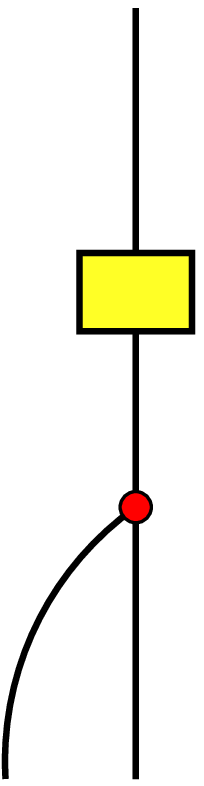}}}
  \put(250,0)     {\scalebox{.36}{\includegraphics{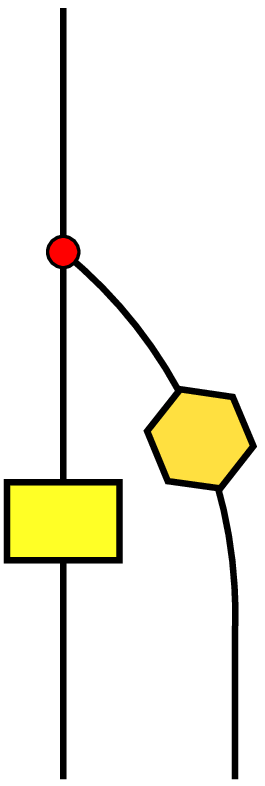}}}
  \put(342,0)     {\scalebox{.36}{\includegraphics{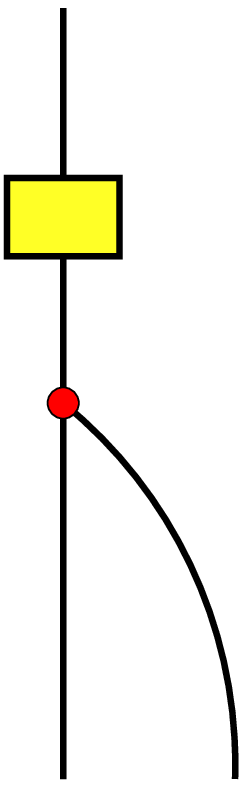}}}
  \put(14.8,32.8) {\small$ \varphi $}
  \put(51.5,42)   {$ = $}
  \put(98.5,63.4) {\small$ \varphi $}
  \put(168,42)    {and}
  \put(253,33.1)  {\small$ \varphi $}
  \put(268.2,42)  {\begin{turn}{40}\small$ \psi $\end{turn}}
  \put(306,50.4)  {$ = $}
  \put(345.6,73.2){\small$ \varphi $}
\eP\\[.3em]
Composing the second of these equalities with 
$(\varphi^{-1}{\circ}\,\eta)\oti\idA$ and using the first of the
equalities twice, one shows that
\\
\begin{picture}(30,88)(-2.5,0)
\setulen79
  \put(0,0)       {\scalebox{.36}{\includegraphics{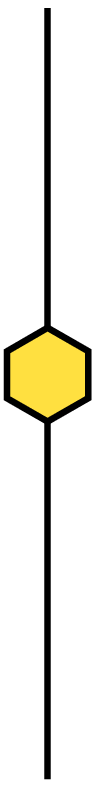}}}
  \put(73.5,0)    {\scalebox{.36}{\includegraphics{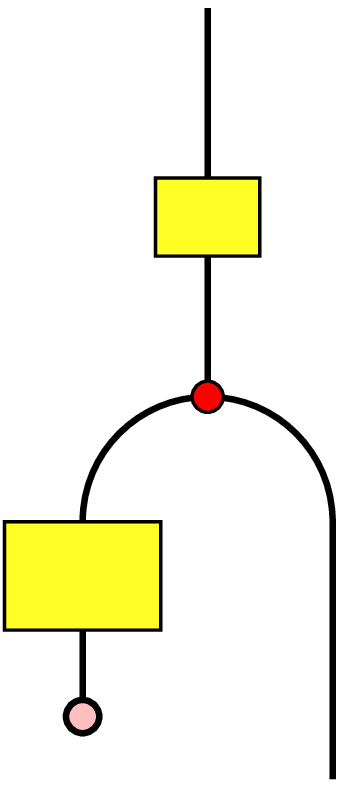}}}
  \put(177.5,0)   {\scalebox{.36}{\includegraphics{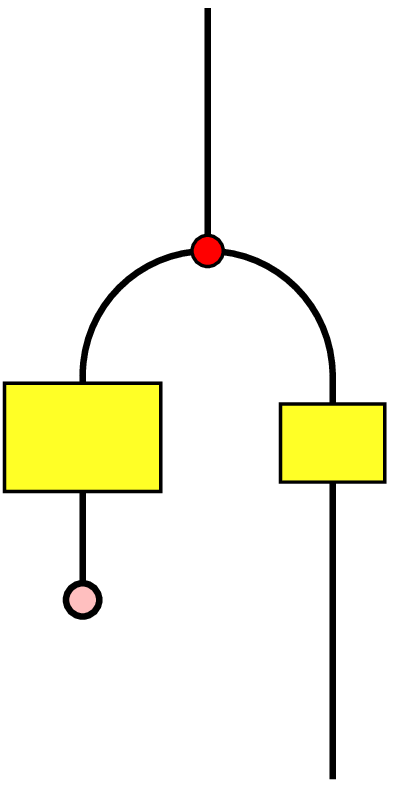}}}
  \put(282,0)     {\scalebox{.36}{\includegraphics{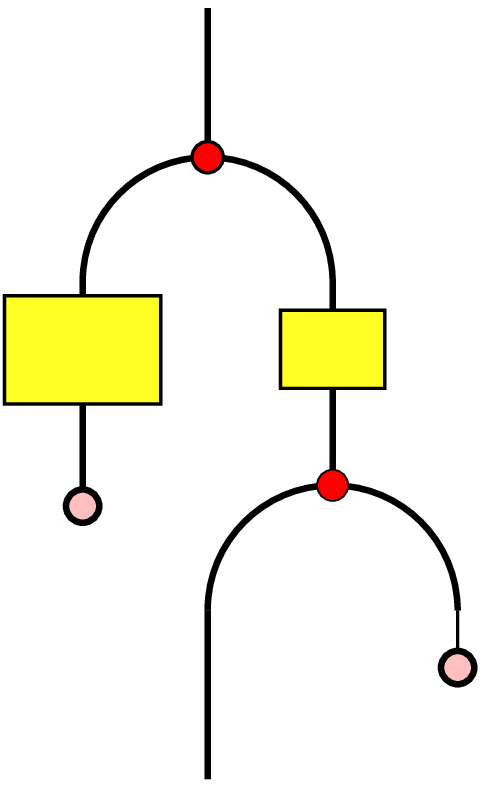}}}
  \put(399,0)     {\scalebox{.36}{\includegraphics{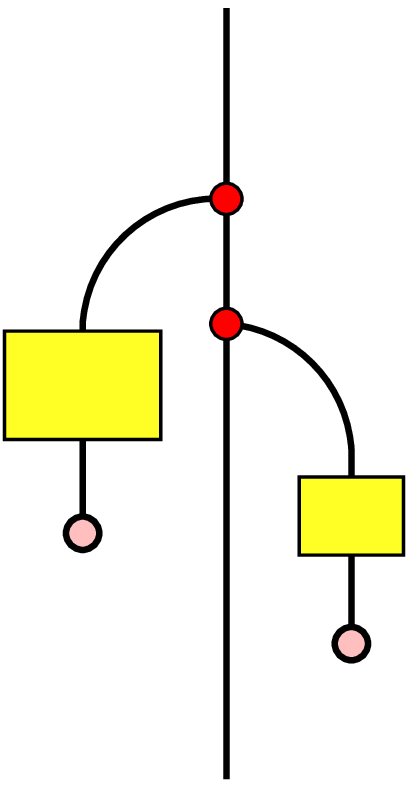}}}
  \put(2.4,51.1)  {\small$ \psi $}
  \put(38.8,42)   {$ = $}
  \put(75.6,24.3) {\small$ \varphi^{\!-1}_{} $}
  \put(95.8,73.4) {\small$ \varphi $}
  \put(141.1,42)  {$ = $}
  \put(179.5,42.4){\small$ \varphi^{\!-1}_{} $}
  \put(216.8,43.5){\small$ \varphi $}
  \put(250,42)    {$ = $}
  \put(283.6,54.4){\small$ \varphi^{\!-1}_{} $}
  \put(321.6,55.2){\small$ \varphi $}
  \put(366,42)    {$ = $}
  \put(400.7,49.1){\small$ \varphi^{\!-1}_{} $}
  \put(440.2,33.8){\small$ \varphi $}
\eP\\[.1em]
Thus $\psi\eq\inn{\varphi^{-1}{\circ}\eta}$; in particular $\psi$ is inner.
\end{proof}

Because of part (i) of the lemma we may restrict our attention to the particular
bimodules $\aAa\idss\psi$. This yields the following result, which is a variant 
of proposition 3.14 of \cite{vazh} (the latter is restated as corollary 
\ref{cor.RZ} below):

\begin{proposition}
\label{prop:ex.seq}
There is an exact sequence
  \be
  1 \,\to\, \Inn(A) \,\to\  \Aut(A)  \ \to\, \PicCAA
  \label{ex.seq}\ee
of groups.  
\end{proposition}

\begin{proof}
{}From (i) and (ii) of lemma \ref{lem.aAa} we deduce that 
$\aAa\idss\psi \otA \aAa\idss{\psi'} \,{\cong}\, \aAa\idss{\psi\psi'}$. Thus the
mapping $\psi\,{ \mapsto}\, [\aAa\idss\psi]$ from $\Aut(A)$ to the Picard group 
$\PicCAA$ of the \cat\ of $A$-bimodules is a group homomorphism. According to 
part (iv) of the lemma, the kernel of this homomorphism is $\Inn(A)$.  
\end{proof}

Note that, up to isomorphism, $\PicCAA$ depends only on the Morita class of
$A$, whereas different representatives of a Morita class can have rather
different \auto\ groups. 
Also, the last morphism in \erf{ex.seq} is not, in 
general, an epimorphism. In the CFT context this means that in general not 
every internal symmetry can be detected as an \auto\ of the algebra $A$ used 
in the TFT construction; an example of this phenomenon has been presented in 
\cite{ffrs3}.


\section{Azumaya \alg s}

For $A$ an \alg\ in a braided monoidal category \C, the monoidal functors
  $$
  \alpha_\AA^\pm:\quad \C \to \CAA 
  $$
of $\alpha$-{\em induction\/} are defined by
$\alpha_\AA^\pm(f) \,{:=}\, \id_\AA \oti f\in\Hom(A\Oti U,A\Oti V)$ for 
morphisms $f\iN\Hom(U,V)$ and $\alpha_\AA^\pm(V)\,{:=}\, (A\Oti V,m\Oti\id_V,
\rho^\pm)$ for objects $V\iN\obj(\C)$. Here the \rep\ morphisms
$\rho^\pm\,{\equiv}\,\rho_V^\pm\iN\Hom(A\Oti V\Oti A,A\Oti V)$ are given by
  $$
  \rho^+ := (m\oti\id_V) \cir (\id_\AA\oti c_{V,A}) ~\quad{\rm and}
  \quad~ \rho^- := (m\oti\id_V) \cir (\id_\AA\oti(c_{A,V})^{-1})\,,
  $$
respectively. These functors first appeared in the context of subfactors
in \cite{lore} and for symmetric monoidal \cats\ in \cite{pare13} (see also
e.g.\ \cite{boek3,ostr,ffrs}). They are used for introducing the following
concept \cite{pare13,vazh}.

\begin{definition}
\label{def.Azumaya}
$A$ is called an {\em Azumaya \alg\/} in \C\ iff the functors $\alinpm$
are monoidal equivalences.
\end{definition}

By remark 3.9 of \cite{ffrs}, the left and right centers of an Azumaya \alg\ 
are isomorphic to $\one$. Also, one has $\alinpm(\one)\eq A$, implying that 
every Azumaya \alg\ is simple. Further, for Azumaya \alg s the functors 
$\alinpm$ induce a group isomorphism (in fact, in general even two different 
isomorphisms) $\PicC \,{\stackrel\cong\longrightarrow}\, \PicCAA$, so that 
from proposition \ref{prop:ex.seq} we obtain the Rosenberg\hy Zelinsky type 
sequence that was established in proposition 3.14 of \cite{vazh}:

\begin{corollary}
\label{cor.RZ}
If $A$ is Azumaya, then there is an exact sequence
  $$ 
  1 \,\to\, \Inn(A) \,\to\  \Aut(A)  \ \to\, \PicC \,.
  $$ 
\end{corollary}

\noindent
Note that for $\C\eq\mathcal Vect_\koerper$, this reduces to Noether's classical
result that all \auto s of a central simple \koerper-algebra are inner.
For non-Azumaya \alg s, $\alinpm$ still induce group homomorphisms from
$\PicC$ to $\PicCAA$, but in general they are neither injective nor surjective
(counter examples from RCFT are given by the \mtcs\ of the D$_{\rm even}$ series 
of $\mathfrak{sl}(2)$ Wess\hy Zu\-mi\-no\hy Wit\-ten models).

\medskip

For \alg s in a braided monoidal \cat\ one defines the {\em opposite\/}
\alg\ of $A$ as $(A,m \cir c_{A,A}^{-1}, \eta)$, and the {\em product\/} of
$A$ and $B$ as \alg s by $(A \Oti B, (m_A \Oti m_B) \cir
  $\linebreak[0]$
(\id_A \oti c_{A,B}^{-1}{\otimes}\,\id_B), \eta_A\Oti\eta_B)$.
Analogously, when both $A\eq(A,m_A,\eta_A,\Delta_A,\eps_A)$ and 
$B\eq(B,m_B,\eta_B,\Delta_B,\eps_B)$ are \ssFA s, one sets
  \be
  \begin{array}l
  \Aop := (\,A,m \cir c_{A,A}^{-1}, \eta, c_{A,A} \cir \Delta, \eps\,) 
  \qquad{\rm and} \\[.3em]
  A \ati B := (\,A\Oti B,(m_A\oti m_B) \cir(\id_A \oti c_{A,B}^{-1}\oti\id_B),
  \eta_A\oti\eta_B, \\ \mbox{\hspace{7.7em}}
  (\id_A\oti c_{A,B}\oti\id_B)\cir(\Delta_A\oti\Delta_B),\,\eps_A\oti\eps_B\,)
  \eear
  \label{oppprod}\ee
and shows (propositions 3.18(ii) and 3.22(ii) of \cite{fuRs4}) that the 
opposite \alg\ \Aop\ and the product \alg\ $A\atimes B$ are again symmetric 
special Frobenius.

If $A$ and $B$ are Azumaya, then so are \Aop\ and $A\atimes B$, and
forming the opposite \alg\ and the product is compatible with Morita equivalence
of \alg s in \C, see theorems 3.3 and 3.4 of \cite{vazh}. Thus the set of Morita
classes of Azumaya \alg s in \C\ carries a natural group structure, with product 
induced by $\atimes$ and inverse induced by taking the opposite \alg; this group
is called the Brauer group of \C.\,%
  \footnote{~%
  In \cite{pare13,vazh} two different Brauer groups are considered.
  They coincide when the tensor unit is projective, and hence in particular
  when \C\ is semisimple.}

\medskip

For the rest of this section, \C\ is again a \mtc. Recall from proposition 
\ref{prop.pertinent} that then the property of an algebra in \C\ of being
simple symmetric special Frobenius is Morita stable. It thus makes sense to 
study the group of Morita classes of symmetric special Frobenius Azumaya 
\alg s in \C; we call this group the {\em \FBg\/} of \C\ and denote it by 
$\Br(\C)$. In the braided setting, $\Br(\C)$ is not necessarily abelian.

The formulas \erf{oppprod} still make sense and are compatible with Morita
equivalence when $A$ (and/or $B$) is not Azumaya. But $A\oti\Aop$ is no longer
Morita equivalent to $\one$ if $A$ is not Azumaya. Thus the set of Morita 
classes of arbitrary \ssFA s in \C\ is no longer a group, though still a 
monoid, see \cite[remark\,5.4]{fuRs4} and \cite{evpi2}.
  
Using notation 
from section 5, we denote by $\Zu(A)$ the square matrix with entries
  $$ 
  \Zu(A)_{ij} := \dimc\,\HomAA(U_i\,\Otip A\,\Otim U_j,A)
  = \dimc\,\HomAA(\alinp(U_i),\alinm(U_j^\vee))
  $$ 
for $i,j\iN\I$. Here the second equality follows by comparing the 
definition of $\otimes^{\!\pm}$ with that of $\alpha$-induction.
One has $\Zu(\one)_{ij}\eq\delta_{i,\j^\vee}$ with
$U_{\j^\vee}\,{\cong}\,U_j^\vee$ as well as (by proposition 5.3 of
\cite{fuRs4}) $\Zu(\Aop) \eq \Zu(A)^{\rm t}$ and
$\Zu(A{\atimes} B) \eq \Zu(A)\,\Zu(\one)\,\Zu(B)$.  

\begin{proposition}
A \ssFA\ $A$ is Azumaya iff $\,\Zu(A)$ is a permutation matrix.
\end{proposition}

\begin{proof}
According to proposition 2.36 of \cite{ffrs} one has
  $$ 
  \begin{array}l
  \HomAA(\alinp(U_i),\alinp(U_j)) = \Hom(U_i,C_l(A)\Oti U_j) \cong
  \complex^{\sum_k\! \N jki \Zu(A)_{k0}}_{\phantom|}\ \quad{\rm and} \\[.3em]
  \HomAA(\alinm(U_i),\alinm(U_j)) = \Hom(U_i,C_r(A)\Oti U_j) \cong 
  \complex^{\sum_k\! \N jki \Zu(A)_{0k}}_{\phantom|} ,  \eear
  $$ 
where $\N jki\eq\dimc\,\Hom(U_j\Oti U_k,U_i)$ and where the equality sign 
indicates a canonical isomorphism (given in (2.69) of \cite{ffrs}), while the 
second isomorphisms follow from remark 3.7 and lemma 3.13 of \cite{ffrs}. Now 
if $\,\Zu(A)\,$ is a permutation matrix, then the number of isomorphism classes
of simple $A$-bimodules equals the number $|\I|$ of isomorphism classes of 
simple objects of \C\ (see remark 5.19(ii) of \cite{fuRs4}). Moreover, because 
of $\alinpm(\one)\eq A$ we have $\Zu(A)_{i0}\eq\delta_{i,0}$, from which it 
follows that $\HomAA(\alinpm(U_i),\alinpm(U_j))\,{\cong}\,\delta_{i,j}\complex$. 
Thus the functors $\alinpm$ are essentially surjective and fully faithful, and 
hence are equivalences.
\\
Conversely, if $\alinpm$ are monoidal equivalences, then the bimodules 
$\alinp(U_i)$ and $\alinm(U_i)$ are simple for all $i\iN\I$, and moreover, 
each simple $A$-bimodule is isomorphic to one of the $\alinp(U_i)$ and to 
one of the $\alinm(U_j)$. Thus $\Zu(A)$ is a permutation matrix. 
\end{proof}

\begin{proposition}
\label{prop.FBr}
The \FBg\ $\Br(\C)$ of a \mtc\ \C\ is finite.
\end{proposition}

\begin{proof}
Every simple algebra in a \mtc\ \C\ is Morita equivalent to an algebra 
$A$ satisfying $\Hom(\one,A)\eq\complex\eta$ \cite{ostr}. 
And the number of \ssFA s in \C\ having this property is finite 
(proposition 3.6 of \cite{fuRs9}).
\end{proof}

\begin{remark}
In RCFT the integers $\Zu(A)_{ij}$ possess the interpretation of the coefficients 
of the torus partition function in a standard basis, called the basis of 
characters, of one-point conformal blocks on the double of the torus. $\Zu(A)$
can thus be described in terms of maximal extensions of the chiral symmetries
for left and right movers. As an additional information, $\Zu(A)$ encodes an
isomorphism of the fusion rules of the two extensions \cite{mose2}.
The procedure of `extending the chiral symmetries' corresponds to passing to the
modular tensor categories $\C^{{\ell{\rm oc}}}_{C_{l/r}(A)}$ of local modules
\cite{pare23,ffrs} of the left and right centers $C_{l/r}(A)$, respectively.
These two categories are monoidally equivalent (theorem 5.20 of \cite{ffrs}),
so that in particular they have indeed isomorphic Grothendieck rings,
$K_0(\C^{{\ell{\rm oc}}}_{C_l(A)})\,{\cong}\,K_0(\C^{{\ell{\rm oc}}}_{C_r(A)})$.
Moreover, one can lift the algebra $A$ to Azumaya algebras in in 
$\C^{{\ell{\rm oc}}}_{C_{l/r}(A)}$ (see proposition 4.14 of \cite{ffrs}).
In this sense, the two extensions are maximal, and the isomorphism of the
fusion rules is encoded in the `Azumaya part' of $A$.
\\
In particular, when $A$ is Azumaya, then $[A]\,{\mapsto}\,\Zu(A)\Zu(\one)$
is a group homomorphism from $\Br(\C)$ to $\Aut(K_0(\C))$.
\end{remark}


\section{Reversions}

Let us now return to the question of finding a suitable equivalence 
relation for Jandl algebras, such that different algebras in the same class
yield, via the TFT construction, one and the same full CFT, i.e.\ full CFTs 
that, up to possibly certain scalar factors (for details see \cite{ffrs5})
have the same correlators.

We recall from definition \ref{jandl} that a Jandl algebra 
$A\,{\equiv}\,(A,\sigma)$ in a ribbon \cat\ \C\ consists of a \ssFA\ $A$ in \C\
and a {\em reversion\/} $\sigma \iN \Hom(A,A)$, obeying the relations \erf{sig}.
Below we will also need the fact that in a ribbon category the bidual functor
$\mbox{\sl?}^{\vee\vee}$ is naturally equivalent to $\id_\C$, via the 
isomorphisms $\,\delta_U \eq ( \id_{U^{\vee\vee}}^{} \oti d_U ) \cir
[ (c^{}_{U^\vee,U^{\vee\vee}} \cir b^{}_{U^\vee}) \oti \theta_U^{} ]
\iN\Hom(U,U^{\vee\vee})\,$.

On the set of isomorphism classes of Jandl algebras in \C\ an equivalence 
relation can be obtained as follows. Let $(A,\sigma)$ be a simple Jandl algebra,
$M\eq(\dot M,\rho)$ a left $A$-module and $B\iN[A]$ the algebra 
$B\eq M^\vee \!\otA M$. Here $M^\vee$ is the right $A$-module 
$M^\vee\eq(\dot M^\vee,(\id_{\dot M^\vee}\Oti d_A)\cir(\rho^\vee\Oti\idA))$. Via
the reversion one can endow the object $\dot M^\vee$ also with a {\em left\/} 
$A$-module structure, denoted by $M^\sigma$ (definition 2.6 of \cite{fuRs8}):
$M^\sigma \eq \big( \dot M^\vee,(\tilde d_A\Oti\id_{\dot M^\vee})\cir
(\sigma\oti{}^\vee\!(\rho\,{\circ}\,c_{\dot M,A})) \big)$.
Now suppose that there is an isomorphism $g \iN \HomA(M,M^\sigma)$ such that
  \be
  g = \Epsilon\, g^\vee{\circ}\,\delta_{\dot M} 
  \label{ggvee}\ee
for some $\Epsilon\iN\complex$.
Then as shown in proposition 2.16 (combined with theorem 2.14) of \cite{fuRs8},
one can use $g$ to construct a reversion 
  \be
  \sigma_g := r \circ (g\oti g^{-1}) \circ
  (\theta_{\dot M}\oti\id_{\dot M^\vee}) \circ c_{\dot M^\vee\!,\dot M} \circ e 
  \label{sigmag}\ee 
of the algebra $B$ (here $e\iN\Hom(B,\dot M^\vee\Oti\dot M)$ and 
$r\iN\Hom(\dot M^\vee\Oti\dot M,B)$
are the embedding and restriction morphisms for $B$ as a retract of
$\dot M^\vee\Oti\dot M$; see formula (2.60) of \cite{fuRs8}).

\begin{definition}
{\rm (i)}~\,Two Jandl algebras $(A,\sigma)$ and $(B,\tau)$ are called
isomorphic, denoted by $(A,\sigma)\,{\cong}\,(B,\tau)$, iff there 
exists an isomorphism $\varphi{:}\ A\To B$ of $A$ and $B$ as
\ssfa s such that $\varphi\cir\sigma\eq\tau\cir\varphi$.
\\[.2em]
{\rm (ii)}~\,Two Jandl algebras $(A,\sigma)$ and $(B,\tau)$ are called
equivalent, denoted by
  \be
  (A,\sigma) \,\simj\, (B,\tau) \,,
  \label{simj}\ee
iff $(B,\tau)$ is isomorphic to $(M^\vee \!\otA M,\sigma_g)$ with $M$ and
$\sigma_g$ as above.
\end{definition}


\begin{lemma}
The relation \,`\,$\simj$'\, {\rm\erf{simj}} is an equivalence relation.
\end{lemma}

\begin{proof}
That \erf{simj} is symmetric is shown by proposition 2.17 \cite{fuRs8}.
\\
Reflexivity is seen by taking $M \eq A$ and $g\eq\Phi\cir\sigma^{-1}$ with
$\Phi$ the isomorphism given in \erf{Phi}. Then e.g.\ the equality \erf{ggvee}
(with $\Epsilon\eq1$) follows with the help of $\Phi\eq\Phi'$ (which
holds because $A$ is symmetric) and of $\eps\cir m\cir(\idA\oti\sigma) \eq
\eps\cir m\cir(\sigma\oti\idA)$.
\\
Transitivity is shown by verifying (which is a bit lengthy, but straightforward) 
that when $(A,\sigma)\,{\simj}\,(B,\tau)$ via an $A$-$B$-bimodule $M$ and 
isomorphism $g$, and $(B,\tau)\,{\simj}\,(C,\varrho)$ via a $B$-$C$-bimodule 
$N$ and isomorphism $h$, then $(A,\sigma)\,{\simj}\,(C,\varrho)$ via the
$A$-$C$-bi\-mo\-dule $M\,{\otimes_B}\,N$ and the morphism 
$c_{\dot M^\vee\!,\dot N^\vee}{\circ}\, (g \oti h)$, which is an isomorphism
between $M\,{\otimes_B}\,N$ and $(M\,{\otimes_B}\,N)^\sigma$.
\end{proof}

It can be checked that the equivalence relation \erf{simj} indeed refines 
Morita equivalence. In particular, for the allowed interpolating bimodules $M$ 
the functor $F_M \eq M{\otimes_B^{}}{-}\,{:}\ \CB\To\CA$ has the property that 
$F_M\cir F_\tau$ and $F_\sigma\cir F_M$ are naturally isomorphic, 
where $F_\sigma$ is the endofunctor of \CA\ that acts on objects
as $X\,{\mapsto}\,X^\sigma$ and on morphisms (when regarded as morphisms of \C)
as $f\,{\mapsto}\,f^\vee$, and analogously for the endofunctor $F_\tau$ of \CB.

Let us also remark that the class of objects $\dot M$ for which there can exist
module isomorphisms satisfying \erf{ggvee} is quite restricted. As shown in
\cite{fuRs8}, the number $\Epsilon$ appearing in \erf{ggvee} can only take the
values $\pm1$. Moreover, by comparison with formula (3.29) of \cite{fuSc16}
one learns that every simple subobject of $\dot M$ has $\Epsilon$ as the value
of its \fsi, i.e.\ this value must be the same for all subobjects 
(also, $\Epsilon$ does not depend on $g$).

Next we note that the existence of reversions does not behave nicely with 
respect to forming the product of algebras. It is easy to check that for 
given Jandl algebras $(A,\sigma)$ and ($B,\sigma')$, the morphism 
$\sigma\oti\sigma'$ is not, in general, a reversion of the product algebra 
$A\ati B$. In fact, $A\ati B$ may not possess any reversion at 
all: owing to $\Zu(\Aop) \eq \Zu(A)^{\rm t}$, existence of a reversion implies 
that the matrix $\Zu(A)$ is symmetric, and this property is not preserved under 
$\atimes$. Consider for instance the case that $A$ and $B$ are Azumaya and 
Jandl. Then $\Zu(A)_{ij}\eq \delta_{i^\vee\!,\pi_A(j)}$ with $\pi_A$ a 
permutation of order 2, and analogously for $B$. But $\Zu(A{\atimes}B)$ is 
given by the permutation $\pi_{\!A\atimes B}\eq\pi_A\cir\pi_B$,
which has order 2 only if $\pi_A$ and $\pi_B$ commute.

As a consequence, an equivalence relation (such as \erf{simj}) between
Jandl algebras is not (in general) compatible with $\atimes\,$.
In particular, $\atimes$ does not induce a group structure on the set of 
equivalence classes of Azumaya\hy Jandl algebras. Thus there appears to be no 
braided analogue of the {\em involutive Brauer group\/} that can be defined 
for algebras \cite{pasr} or for symmetric monoidal \cats\ \cite{vevi2}.

\medskip

As already pointed out, Morita equivalent \ssfa s yield, via the TFT 
construction, the same {\em oriented\/} full CFT. Comparing the correlators of 
the {\em unoriented\/} full CFTs that are obtained from different Jandl 
algebras is a much more difficult task than in the oriented case, and we have 
not yet fully investigated this issue. But we plan to establish in a subsequent 
paper \cite{ffrs5} that it is indeed the equivalence relation \erf{simj} that 
takes over the role of Morita equivalence in unoriented CFT. There is, in fact, 
a special situation in which the comparison of the CFTs obtained from two Jandl 
\alg s $(A,\sigma)$ and $(B,\tau)$ is much simplified, namely when $B\eq A$. 
In this case one can use the fact that for any two reversions $\sigma$ and 
$\sigma'$ of $A$, $\omega\,{:=}\,\sigma^{-1}{\circ}\,\sigma'$ is an algebra
\auto\ of $A$. Conversely, while the composition of a reversion and an algebra
automorphism is not, in general, again a reversion, we still have the following 

\begin{proposition}
\label{prop.oso=s}
Let $(A,\sigma)$ be a Jandl \alg\ and $\omega$ an algebra automorphism of $A$. 
\\[.2em]
{\rm (i)} $(A,\sigma\cir\omega)$ is a Jandl \alg\ iff 
  \be
  \omega\cir\sigma\cir\omega = \sigma \,.
  \label{oso=s}\ee
{\rm (ii)} If $A$ is simple, $\omega\eq\inn\alpha$ is inner, and 
$\inn\alpha\cir\sigma\cir \inn\alpha \eq \sigma$, then
  \be
  (A,\sigma) \,\simj\, (A,\sigma\cir\inn\alpha) \,.
  \ee
\end{proposition}

\begin{proof}
(i) is shown in proposition 2.3 of \cite{fuRs8}.
\\[.2em]
(ii) Consider the morphism 
  $$
  g_\alpha := \epsa\,\Phi\cir\sigma^{-1}{\circ}\,m\cir(\idA\oti\alpha^{-1}) \
  \in\Hom(A,A^\vee) \,, 
  $$
with $\Phi$ as in \erf{Phi} and $\epsa\iN\{\pm1\}$.
By using the defining properties of $A$ (being a \ssfa) and of $\omega$
and $\sigma$, it is not difficult to see that $g_\alpha$ intertwines the
left $A$-modules $A$ and $A^\sigma$. Similarly, making in addition use of the 
equality \erf{oso=s} one shows that $g_\alpha$ satisfies \erf{ggvee} (with 
$\Epsilon\eq1$). Finally, regarding $A\,{\cong}\,A^\vee\Ota A$ as a retract of 
$A^\vee\Oti A$ via $e\,{:=}\,(\id_{A^\vee}\oti m)\cir(\tilde b_A\oti\idA)$ and 
$r\,{:=}\,(d\oti\idA)\cir(\id_{A^\vee}\oti\Delta)$, one can show that the 
morphism in \erf{sigmag} becomes $\sigma_{g_\alpha} \eq \sigma\cir\inn\alpha$,
provided that
  $$
  \sigma\cir\alpha = \epsa\,\alpha \,.
  $$
As will be seen in proposition \ref{prop.19} below, this equality is 
equivalent to $\inn\alpha\cir\sigma\cir \inn\alpha \eq \sigma$. 
\\[2pt]
Thus indeed $(A,\sigma) \,{\simj}\, (A,\sigma\cir\inn\alpha)$, via $M\eq A$ and
$g\eq g_\alpha$.
\end{proof}

According to proposition \ref{prop.oso=s}, reversions of a \ssFA\ $A$ that 
endow $A$ with inequivalent Jandl structures are related by {\em outer\/} 
algebra \auto s of $A$. Thus the number of inequivalent Jandl structures on
$A$ is bounded by the order of the group $\Aut(A)/\Inn(A)$ of outer \auto s,
which in turn by proposition \ref{prop:ex.seq} is bounded by the order of
\PicCAA. Now for any \ssFA\ in a \mtc\ \C, $\PicCAA$ is a finite group. 
When combined with (the proof of) proposition \ref{prop.FBr} and with the 
results of section 2.5 of \cite{fuRs8}, this yields

\begin{corollary}
The number of equivalence classes, with respect to the relation {\rm\erf{simj}}, 
of simple Jandl algebras in a \mtc\ is finite.
\end{corollary}

In accordance with our remarks about the equivalence of unoriented CFTs above,
part (ii) of proposition \ref{prop.oso=s} leads to the following counterpart 
in CFT, which will be demonstrated in \cite{ffrs5}:

\begin{proposition}
Let $\inn\alpha$ be an inner \auto\ of a simple \ssfa\ $A$ in a \mtc\ \C\ and
$\sigma$ a reversion of $A$ such that $\sigma\cir\inn\alpha$ is a reversion, 
too. Then the correlation functions of the full CFTs based on the Jandl algebras
$(A,\sigma)$ and $(A,\sigma\cir\inn\alpha)$ differ at most by a sign:
  $$
  \Cf_{(A,\sigma)}(\X)
  = (\epsa)^{cr(\X)}_{}\,\Cf_{(A,\sigma\circ\inn\alpha)}(\X) \,.
  $$
$cr(\X)$ is the `number of crosscaps' of the world sheet \X,
which is defined modulo \mbox{\sl2}.
\end{proposition}

The TFT construction of the correlators of unoriented full CFTs involves, as 
compared to the oriented case, in addition at each vertex of the triangulation 
of \iox\ a choice of local orientation. Further, the reversion $\sigma$
enters the prescription only via those edges of the triangulation for which the
chosen orientation at the vertex on one end, when transported along the
edge to its other end, is different from the orientation chosen for the vertex 
on that end. The proof then relies on the fact that the local orientations can 
be chosen in such a way that there are at most two such particular edges
(e.g.\ none of them if \X\ is orientable). This way the effects of replacing
$\sigma$ by $\inn\alpha\cir\sigma$ are confined to neighborhoods of at
most two $A$-ribbons in \iox, and can therefore be analyzed relatively easily. 

\medskip

For completing the proof of proposition \ref{prop.oso=s}(ii), we must still 
study the question of when together with $\sigma$ also $\sigma\cir\inn\alpha$ 
is a reversion. 

\begin{lemma}
\label{lem.rev1}
Given a reversion $\sigma$ of a \ssfa\ $A$ and an invertible $\alpha\iN\Atop$,
the morphism $\sigma\cir\inn\alpha$ is again a reversion of $A$ iff
$\inn{\sigma\circ\alpha} \eq \inn\alpha$.
\end{lemma}

\begin{proof}
According to proposition \ref{prop.oso=s}(i), $\sigma\cir\inn\alpha$ is a 
reversion iff the equality $\sigma\cir\inn\alpha\cir\sigma^{-1}\eq \inn\alpha
^{-1}\eq\inn{\alpha^{-1}}^{}$ holds. On the other hand, by using the defining 
properties \erf{sig} of a reversion one shows $\sigma\cir\inn\alpha\cir\sigma
^{-1}\eq \inn{\sigma\circ\alpha^{-1}}$. Thus $\inn{\alpha^{-1}}\eq\inn{\sigma
\circ\alpha^{-1}}$, which is equivalent to $\inn{\sigma\circ\alpha}\eq\inn\alpha$.
\end{proof}

Note that, using lemma \ref{lem.div.inn}, $\inn{\sigma\circ\alpha}\eq\inn\alpha$
iff $\alpha^{-1}{\multo}\, (\sigma{\circ}\alpha) \iN \Ctopx$. Also, for any 
\ssFA\ $A$ one has $\Hom(\one,C_l(A)) \,{\cong}\, \HomAA(A,A)$ (proposition 2.36
of \cite{ffrs}). Now assume that $A$ is simple. Then $\HomAA(A,A)\,{\cong}\,
\koerper$, and
it follows that $\Ctopx\eq\koerper^{\times}\eta$. Hence 
$\inn{\sigma\circ\alpha}\eq\inn\alpha$ means that $\sigma\cir\alpha$ is a 
nonzero multiple of $\alpha$. Furthermore, owing to 
$\sigma\cir\sigma\cir\alpha\eq\theta_\AA\cir\alpha\eq\alpha$, this multiple 
must be $\pm1$. Further, because of $\Ctopx\eq\koerper^{\times}\eta$ we also
learn that, for simple $A$, $\inn\alpha\eq\inn{\alpha'}$ iff $\alpha\eq\zeta
\alpha'$ for some $\zeta\iN\koerper^{\times}$.
Thus we have proven

\begin{proposition}
\label{prop.19}
Let $(A,\sigma)$ be a simple Jandl \alg\ and $\alpha\iN\Atopx$.
\\[.2em]
{\rm(i)}\, The endomorphism $\sigma\cir\inn\alpha$ is a reversion of $A$ iff 
$\,\sigma\cir\alpha \eq \epsa\,\alpha$ with $\epsa\iN\{\pm1\}$.
\\[.2em]
{\rm(ii)} If $\inn\alpha\eq\inn{\alpha'}$, then $\epsa\eq\epsap$.
\end{proposition}

\newpage 

 \newcommand\wb{\,\linebreak[0]} \def\wB {$\,$\wb}
 \newcommand\Bi[2]    {\bibitem[#2]{#1}}
 \newcommand\JO[6]    {{\em #6}, {#1} {#2} ({#3}), {#4--#5} }
 \newcommand\J[7]     {{\em #7}, {#1} {#2} ({#3}), {#4--#5} {{\tt [#6]}}}
 \newcommand\JJ[6]    {{\em #6}, {#1} {#2} ({#3}), {#4} {{\tt [#5]}}}
 \newcommand\Pret[2]  {{\em #2}, pre\-print {\tt #1}}
 \newcommand\BOOK[4]  {{\em #1\/} ({#2}, {#3} {#4})}
 \newcommand\inBO[8]{{\em #8}, in:\ {\em #1}, {#2}\ ({#3}, {#4} {#5}), p.\ {#6--#7}}
 \def\dim   {dimension}
 \def\jf    {J.\ Fuchs}
 \def\adma  {Adv.\wb Math.}
 \def\coma  {Con\-temp.\wb Math.}
 \def\Coma  {Con\-temp. Math.}
 \def\comp  {Com\-mun.\wb Math.\wb Phys.}
 \def\Comp  {Com\-mun.\wb Math. Phys.}
 \def\cpma  {Com\-pos.\wb Math.}
 \def\duke  {Duke\wB Math.\wb J.}
 \def\fiic  {Fields\wB Institute\wB Commun.}
 \def\Fiic  {Fields Institute\wB Commun.}
 \def\fiiC  {Fields\wB Institute Commun.}
 \def\ihes  {Publ.\wb Math.\wB I.H.E.S.}   
 \def\inma  {Invent.\wb math.}
 \def\joal  {J.\wB Al\-ge\-bra}
 \def\jomp  {J.\wb Math.\wb Phys.}
 \def\josp  {J.\wb Stat.\wb Phys.}
 \def\jpaa  {J.\wB Pure\wB Appl.\wb Alg.}
 \newcommand\ncmp[3] {\inBO{IXth International Congress on
            Mathematical Physics} {B.\ Simon, A.\ Truman, and I.M.\ Da\-vis, 
            eds.} {Adam Hilger}{Bristol}{1989} {#1}{#2}{#3} } 
 \def\nuci  {Nuovo\wB Cim.}
 \def\nupb  {Nucl.\wb Phys.\ B}
 \def\Nupb  {Nucl. Phys.\ B}
 \def\phrl  {Phys.\wb Rev.\wb Lett.}
 \def\pnas  {Proc.\wb Natl.\wb Acad.\wb Sci.\wb USA}
 \def\pspm  {Proc.\wb Symp.\wB Pure\wB Math.}
 \def\rims  {Publ.\wB RIMS}
 \def\rvmp  {Rev.\wb Math.\wb Phys.}
 \def\slnm  {Sprin\-ger\wB Lecture\wB Notes\wB in\wB Mathematics}
 \newcommand\Slnm[1] {{\rm[\slnm\ #1]}}
 \def\trgr  {Trans\-form. Groups}
   \def\AMS    {{American Mathematical Society}}
   \def\AW     {{Addi\-son\hy Wes\-ley}}
   \def\BIR    {{Birk\-h\"au\-ser}}
   \def\Be     {{Berlin}}
   \def\Bo     {{Boston}}
   \def\Ca     {{Cambridge}}
   \def\CUP    {{Cambridge University Press}}
   \def\MD     {{Marcel Dekker}}
   \def\NY     {{New York}} 
   \def\OUP    {{Oxford University Press}}
   \def\PR     {{Providence}}
   \def\SV     {{Sprin\-ger Ver\-lag}}
   \def\WS     {{World Scientific}}
   \def\Si     {{Singapore}} 

\bibliographystyle{amsalpha} \medskip
\end{document}